\documentclass[12pt]{article}
\usepackage{latexsym,amssymb,amsmath,a4wide,color, graphicx}
\usepackage{stackengine}

\newtheorem{thm}{Theorem}
\newtheorem{lem}{Lemma}
\newtheorem{cor}{Corollary}
\newtheorem{conj}{Conjecture}
\newtheorem{prop}{Proposition}
\newtheorem{exer}{Exercise}

\newcommand{\pr}{{\bf Proof.}\ }
\newcommand{\bt}{\begin{thm}}
\newcommand{\et}{\end{thm}}
\newcommand{\bl}{\begin{lem}}
\newcommand{\el}{\end{lem}}
\newcommand{\bp}{\begin{prop}}
\newcommand{\ep}{\end{prop}}
\newcommand{\bc}{\begin{cor}}
\newcommand{\ec}{\end{cor}}
\newcommand{\bcj}{\begin{conj}}
\newcommand{\ecj}{\end{conj}}
\newcommand{\bex}{\begin{exer}}
\newcommand{\eex}{\end{exer}}
\newcommand{\bi}{\begin{itemize}}
\newcommand{\ei}{\end{itemize}}
\newcommand{\be}{\begin{equation}}
\newcommand{\ee}{\end{equation}}
\newcommand{\ben}{\begin{enumerate}}
\newcommand{\een}{\end{enumerate}}

\newcommand{\mt}{t\kern-0.035cm\char39\kern-0.03cm}
\newcommand{\ml}{l\kern-0.035cm\char39\kern-0.03cm}
\newcommand{\md}{d\kern-0.035cm\char39\kern-0.03cm}

\newcommand\xrowht[2][0]{\addstackgap[.5\dimexpr#2\relax]{\vphantom{#1}}}

\newcommand{\veps}{\varepsilon}

\newcommand{\GF}{{\rm GF}}

\newcommand{\dia}{{\rm dia}}
\newcommand{\off}{{\rm off}}
\newcommand{\PSL}{{\rm PSL}}
\newcommand{\PGL}{{\rm PGL}}

\newcommand{\GL}{{\rm GL}}

\newcommand{\ovl}{\overline}

\begin{document}

\title{\vspace{-3cm} Characters of twisted fractional linear groups}

\author{}
\date{}
\maketitle

\begin{center}
\vspace{-1.3cm}

{\large So\v{n}a Pavl\'ikov\'a} \\
\vspace{1.5mm} {\small Slovak University of Technology, Bratislava, Slovakia}\\

\vspace{5mm}

{\large Jozef \v Sir\'a\v n} \\
\vspace{1.5mm} {\small
Open University, Milton Keynes, U.K., and \\ Slovak University of Technology, Bratislava, Slovakia}

\vspace{4mm}

\end{center}

\begin{abstract}
We determine character tables for twisted fractional linear groups that form the `other' family in Zassenhaus' classification of finite sharply 3-transitive groups.

\end{abstract}

\vskip 3mm

\section{Introduction and preliminaries}\label{sec:pre}

A classification of all finite, sharply 3-transitive permutation groups follows from a classical result by Zassenhaus \cite{Zas}, by which the only two families of such groups are the fractional linear groups $\PGL(2,q)$ for any prime power $q$ and their `twisted' companions $M(q^2)$ for any odd prime power $q$. Recall that $\PGL(2,q)$ can be introduced as the group of fractional transformations $z\mapsto (az+b)/(cz+d)$ of the set $\GF(q)\cup\{\infty\}$, where $ad-bc\ne 0$, with the obvious rules for calculations with $\infty$. In the case of a finite field $F$ of the form $F=\GF(q^2)$ for an odd prime power $q$ one may `twist' the fractional transformations by considering the permutations of $F\cup \{\infty\}$ defined by $z\mapsto (az+b)/(zc+d)$ if $ad-bc\in S(F)$ and $z\mapsto (az^{\sigma} +b)/(cz^{\sigma}+d)$ if $ad-bc\in N(F)$, where $S(F)$ and $N(F)$ are the sets of non-zero squares and non-squares of $F$ and $\sigma$ is the unique involutory (Galois) automorphism of $F$. The collection of all such `untwisted' and `twisted' fractional transformations under composition constitutes the {\em twisted fractional linear group} $M(q^2)$.
\smallskip

The notation $M(q^2)$ comes from the monographs \cite[p. 261]{Pass} and \cite[p. 163]{HB}, and is also used in the textbooks \cite[p. 188]{Rob} and \cite[p. 283]{Rot}; the letter $M$ was introduced in the original article \cite[p. 36]{Zas} as a tribute to Mathieu who discovered the first group in this series (for $q=3$, of degree $10$). A possible alternative notation could be derived from the existence of just three non-trivial $2$-extensions of $\PSL(2,q^2)$ by outer automorphisms -- a diagonal automorphism $\delta$, the Galois automorphism $\sigma$, and their product $\delta\sigma$ -- leading to the groups $\PGL(2,q^2)$, $P\Sigma L(2,q^2)$ and $M(q^2) \cong \PSL(2,q^2)\langle\delta\sigma\rangle$, see \cite{Whi}.
\smallskip

While the fractional linear groups $\PGL(2,q)$ and their subgroups $\PSL(2,q)$ have been studied widely, their twisted versions $M(q^2)$ have received comparatively less attention, and this also applies to characters. A number of resources on representations and characters of fractional linear groups are accessible (either in form of explicit tables or in terms of methods of their derivation, cf. \cite{Dor,JW}). Character tables for the `other' family of sharply $3$-transitive groups, by contrast, do not appear to be available.
\smallskip

Among numerous applications of characters we mention here the ones related to the theory of regular hypermaps on compact surfaces, or, equivalently, finite groups generated by three involutions. In \cite{Adr}, enumeration of regular hypermaps of a given type on fractional linear groups was obtained with the help of the Frobenius' character formula \cite{Fr} for counting tuples of group elements with entries in given conjugacy classes; more applications of this type can be found in \cite{Jon}. There is, however, no corresponding result for regular hypermaps on {\em twisted} fractional linear groups, although enumeration of regular maps (of unspecified type) on these groups can be found in \cite{EHS}.
\smallskip

The purpose of this paper is to calculate character tables for the twisted fractional linear groups $M(q^2)$ for any odd prime power $q$. This will be done in a completely elementary way by means of standard results in representation theory, by lifting (that is, inducing) characters from subgroups of $M(q^2)$ onto the entire group and finding their decomposition into irreducible constituents.
\smallskip

To avoid fractional transformations we will work with a representation of $M(q^2)$ used in \cite{EHS}), which also differs from the one of \cite{Whi} and which will prove more suitable for our purposes. For $F=\GF(q^2)$, $q$ an odd prime power, we let $J= \GL(2,F) \rtimes \langle\sigma\rangle$, with $\langle\sigma \rangle$ identified with the additive group $C_2$ in the obvious way and with multiplication given by $(A,r)(B,s)= (AB^{\sigma^{r}},r+s)$, where $B^{\sigma}$ is obtained from $B$ by applying $\sigma$ to every entry. For each $A\in \GL(2,F)$ let $\iota_A\in C_2=\{0,1\}$ be defined by $\iota_A=0$ if $\det(A)\in S(F)$ and $\iota_A=1$ if $\det(A)\in N(F)$. The `twisted' subgroup $K$ of $J$ of index $2$ is defined by letting $K=\{(A,\iota_A);\ A\in \GL(2,F)\}$ with multiplication as before, that is,  $(A,\iota_A) (B,\iota_B) =(AB^{\sigma^{\iota_{\scalebox{0.35}{$A$}}}},\iota_A+\iota_B)$ for every $A,B\in \GL(2,F)$.
\smallskip

Consider now the subgroup $K_0=\{(A,0);\ A\in \GL(2,F),\ \iota_A=0\}$ of $K$ of index $2$. The centre $L$ of $K_0$ consists of pairs $(D,0)$, where $D\in \GL(F)$ is a scalar matrix; note that $L$ is also normal in both $K$ and $J$. The factor group $G=K/L$ turns out to be isomorphic to $M(q^2)$. We will identify $G$ with $M(q^2)$ throughout. This way, $G$ can be regarded as a subgroup of index $2$ of the group $\ovl{G}=J/L$, and the factor group $H=K_0/L$ can be identified with $\PSL(2,F)$. Elements $(A,\iota_A)L$, that is, cosets $\{(\delta A,\iota_A);\ \delta\in F^*\}$ of the factor group $G=K/L$, will throughout be denoted $[A,\iota_A]$; they will be called {\em untwisted} if $\iota_A=0$ and {\em twisted} if $\iota_A=1$.

\section{Conjugacy classes of $M(q^2)$}\label{sec:conj}

To be in position to consider characters of the twisted group $M(q^2)=G=K/L$ we will need to determine conjugacy classes of $G$. This was done in \cite{EHS} for conjugacy of twisted elements with respect to $\ovl{G}$, and it turns out that the detailed analysis therein furnishes all one needs to determine the conjugacy within $G$, a subgroup of $\ovl{G}$ of index two. We sum up the corresponding results in what follows, using notation and machinery of \cite{EHS}. For any non-zero elements $a,b$ of a field $\tilde F$ we let $\dia(a,b)$ and $\off(a,b)$ be the $2\times 2$ matrices of ${\rm GL}(2,\tilde F)$ with, respectively, the diagonal and off-diagonal entries $a,b$ and with remaining entries equal to zero.
\smallskip

We begin by conjugacy of untwisted elements of $G=M(q^2)$, forming the subgroup $H\cong \PSL(2,q^2)$ of $G$. Every untwisted element of $G$ can be identified with $[A,0]$ for $A\in \PSL(2,q^2)$, or, for short, just with $A$, tacitly assuming that writing $A\in \PSL(2,q^2)$ means $\pm A$ for $A\in {\rm SL}(2,q^2)$ as usual. Since conjugacy classes in $\PSL(2,q^2)\cong H$ are well known, the question is which pairs are fused by conjugacy in $G$. Let $\xi$ be a primitive element of $F=\GF(q^2)$ and let $\zeta$ be a primitive $(q^2+1)^{\rm th}$ root of unity in an extension $F'$ of $F$ of degree two. Further, let $B_1$ and $B_{\veps}$ be the elements of $\PSL(2,q^2)$ obtained from the identity matrix by replacing the top right $0$ with $1$ and with some $\veps\in N(F)$, respectively. Our calculations will, in spirit, be similar to those in \cite[Lemma 4.7]{Whi}.
\smallskip

Now, a non-identity element $A\in \PSL(2,q^2)$ is conjugate within $\PSL(2,q^2)\cong H$ to either $w=\off(1,-1)$, or to one of $u=B_1$ and $u'=B_{\veps}$, or else to $a(\theta) = \dia(\theta,\theta^{-1})$ for some $\theta$ that is a power of $\xi$ or $\zeta$. In the first three cases the elements $w$, $u$ and $u'$ generate a single conjugacy class each. The elements $a(\theta)$ for $\theta$ in $F$ and $F'$ generate, respectively, $(q^2-5)/4$ distinct conjugacy classes for $\theta =\xi^j$, $1\le j\le (q^2-5)/4$, and another $(q^2-1)/4$ distinct conjugacy classes for $\theta = \zeta^j$, $1\le j\le (q^2-1)/4$. The corresponding centralizers in $G$ of the above five types of elements have orders $q^2-1$,  $q^2$, $q^2$, $(q^2-1)/2$ and $(q^2+1)/2$, respectively. (We included this information also in Table \ref{tab:PSL} in Section \ref{sec:pre-char} displaying irreducible characters of $H$.)
\smallskip

With the help of this we determine conjugacy classes of untwisted elements in the overgroup $G$ of $H$ of index two. Conjugacy in $G$ fuses the classes generated by the untwisted elements $[u,0]$ and $[u',0]$ into a single class (e.g. by conjugating by the twisted element $[\dia(\veps,1),1]$), still with centralizer of order $q^2$, while the class generated by $[w,0]$ remains the same under conjugacy in $G$, with centralizer of twice the original order, that is, $2(q^2-1)$. Sorting out the untwisted elements $[a(\theta),0]$ needs more care. Note first that the traces (defined up to multiplication by $-1$) of $A,A'\in {\rm SL}(2,q^2)$ forming two untwisted elements $[A,0]$ and $[A',0]$ which are conjugate by a twisted element satisfy ${\rm tr}(A')= \pm({\rm tr}(A))^{\sigma}$. Observe also that the twisted element $[\dia(\veps,1),1]$ conjugates the untwisted element $[a(\theta),0]$ to $[a(\theta^{\sigma}),0]$, and note that $a(\theta)$ is conjugate to $a(\theta^{-1})$ in $\PSL(2,q^2)$.
\smallskip

This implies that two conjugacy classes in $H$ of untwisted elements $[a(\theta),0]$ and $[a(\theta^{\sigma}),0]$ are fused by conjugation in $G$ unless $[a(\theta^{\sigma}),0]$ itself is contained in the $H$-conjugacy class of $[a(\theta),0]$. But the latter means that $a(\theta)$ and $a(\theta^{\sigma})$ are conjugate elements of $\PSL(2,q^2)$, which is if and only if the traces of the two elements are the same up to a sign, that is, $\theta^{\sigma}+ (\theta^{\sigma})^{-1} = \pm(\theta + \theta^{-1})$. Since taking $\sigma$-images means raising in the power of $q$, this equation reduces to $\theta^q+\theta^{-q}\mp (\theta + \theta^{-1})=0$, which is equivalent to $(\theta^{q+1}\mp 1)(\theta^{q-1}\mp 1)=0$. This happens if and only if $\theta$ is one of the $(q-1)^{\rm th}$ or $(q+1)^{\rm th}$ root of $1$ or $-1$ in $F={\rm GF}(q^2)$; in each case there are, respectively, $q-1$ and $q+1$ such roots in $F$. These are precisely the situations when the $G$-conjugacy class of $[a(\theta),0]$ coincides with that of $[a(\theta^{\sigma}),0]$; in all the remaining cases, that is, when $\theta^{\sigma}=\theta^q\notin \{\pm \theta, \pm\theta^{-1}\}$, the two classes are fused by conjugacy in $G$.
\smallskip

Taking this further, if a representative $[a(\theta),0]$ of a non-fused conjugacy class has been chosen for a particular $\theta$ we may automatically disregard the values $-\theta$ and $\pm \theta^{-1}$. Note that if $\theta^2=-1$ then the matrices $a(\theta)$ and $w$ are conjugate in $\PSL(2,q^2)$, and if $\theta^2=1$ then $a(\theta)$ represents the identity element of $\PSL(2,q^2)$. Excluding the values of $\theta$ for which $\theta^4=1$, it follows that in the case when $q\equiv 1$ mod $4$ this leaves, respectively, only $(q-5)/4$, $(q-1)/4$, $(q-1)/4$, and $(q-1)/4$ values of $\theta$ such that $\theta^{q-1}=1$, $\theta^{q-1}=-1$, $\theta^{q+1}=1$, and $\theta^{q+1}=-1$, with the property that the conjugacy classes of elements $[a(\theta),0]$ considered above are mutually distinct. Similarly, if $q\equiv 3$ mod $4$, there are only $(q-3)/4$, $(q-3)/4$, $(q-3)/4$, and $(q+1)/4$ values of $\theta$ with $\theta^{q-1} =1$, $\theta^{q-1}=-1$, $\theta^{q+1}=1$, and $\theta^{q+1}=-1$, respectively, with distinct conjugacy classes of elements $[a(\theta),0]$ as above. In both cases one has a total of $4q-8$ such {\sl distinct} values of $\theta$ giving distinct conjugacy classes; recall that the four $4$th roots of $1$ do not contribute to the classes considered here.
\smallskip

The set $S=\{\theta \in F;\ \theta ^{q\pm 1}\ne \pm 1\}$ is therefore of size $(q^2{-}1)- (4q{-}8)-4 = (q-1) (q-3)$; observe that this number is a multiple of $8$ as $q^2\equiv 1$ mod $4$. The important property of the set $S$ is that it admits a partition into subsets of size $8$ of the form $\{\pm\theta,\pm \theta^{-1},\pm \theta^q,\pm \theta^{-q} \}$; the fact that all these $8$ elements are distinct follows from the way $S$ has been defined. Let $S'$ be an arbitrary but fixed set of distinct representatives of the partition just described, so that $|S'|=|S|/8= (q-1)(q-3)/8$, and for a fixed primitive element $\xi\in F$ let ${\cal U}=\{j; \ 1{\le}j{\le} q^2{-}2;\ \xi^j \in S'\}$. Similarly, consider the set $T=\{\theta \in F'{\setminus} \{\pm 1\};\ \theta^{q^2+1}=1\}$ with $|T|=q^2-1$, which is again a multiple of $8$ (exclusion of those $\theta$ for which $\theta^2=\pm 1$ now reduces to leaving out $\pm 1$ only); all $\theta\in T$ are powers of the primitive $(q^2+1)^{\rm th}$ root of unity $\zeta$ introduced earlier. The set $T$ also admits a partition into $8$-element subset of exactly the same shape as before, and we let $T'$ be any fixed set of distinct representatives of this partition, with $|T'|=(q^2-1)/8$. Finally, let ${\cal V}=\{j;\ 1{\le}j{\le}q^2{-}2;\ \zeta^j\in T'\}$.
\smallskip

Using the above facts and working out the values of $\theta$ in the case $\theta ^{q\pm 1}=\pm 1$ as powers of $\xi$ we arrive at Table \ref{tab:untwisted} for the $1{+}(q{+}1)^2/4$ conjugacy classes of untwisted elements of $G$.


\begin{table}[ht]
	\centering
\begin{tabular}{|c||c|c|c|c|}
\hline \xrowht[()]{4pt}
Untwisted representatives & \# Classes & $|{\rm Class}|$ & $|C_G|$ & Notes\\ \hline \hline \xrowht[()]{4pt}
$[I,0]$ & $1$ & $1$ & $q^2(q^4{-}1)$ & -- \\ \hline \xrowht[()]{4pt}
$[u,0]$ & $1$ & $q^4-1$ & $q^2$ & -- \\ \hline \xrowht[()]{4pt}
$[w,0]$ & $1$ & $q^2(q^2{+}1)/2$  & $2(q^2-1)$ & -- \\ \hline\hline \xrowht[()]{4pt}
$[a(\xi^{j(q+1)}),0],\ j{\le}(q{-}5)/4$ & $(q-5)/4$ & $q^2(q^2+1)$ & $q^2-1$ & $q\equiv 1\ {\rm mod}\ 4$\\ $[a(\xi^{j(q+1)}),0],\ j{\le}(q{-}3)/4$ & $(q-3)/4$ & $q^2(q^2+1)$ & $q^2-1$ & $q\equiv 3\ {\rm mod}\ 4$\\ \hline \xrowht[()]{4pt}
$[a(\xi^{j(q+1)/2}),0],\, {\rm odd}\ j{\le}(q{-}3)/2$ & $(q-1)/4$ & $q^2(q^2+1)$ & $q^2-1$ & $q\equiv 1\ {\rm mod}\ 4$\\
$[a(\xi^{j(q+1)/2}),0],\, {\rm odd}\ j{\le}(q{-}5)/2$ & $(q-3)/4$ & $q^2(q^2+1)$ & $q^2-1$ & $q\equiv 3\ {\rm mod}\ 4$\\ \hline \xrowht[()]{4pt}
$[a(\xi^{j(q-1)}),0],\ j{\le}(q{-}1)/4$ & $(q-1)/4$ & $q^2(q^2+1)$ & $q^2-1$ & $q\equiv 1\ {\rm mod}\ 4$\\ $[a(\xi^{j(q-1)}),0],\ j{\le}(q{-}3)/4$ & $(q-3)/4$ & $q^2(q^2+1)$ & $q^2-1$ & $q\equiv 3\ {\rm mod}\ 4$\\ \hline \xrowht[()]{4pt}
$[a(\xi^{j(q-1)/2}),0],\, {\rm odd}\ j{\le}(q{-}3)/2$ & $(q-1)/4$ & $q^2(q^2+1)$ & $q^2-1$ & $q\equiv 1\ {\rm mod}\ 4$\\
$[a(\xi^{j(q-1)/2}),0],\, {\rm odd}\ j{\le}(q{-}1)/2$ & $(q+1)/4$ & $q^2(q^2+1)$ & $q^2-1$ & $q\equiv 3\ {\rm mod}\ 4$\\ \hline\hline \xrowht[()]{4pt}
$[a(\xi^j),0],\ j\in {\cal U}$ & $(q{-}1)(q{-}3)/8$ & $2q^2(q^2+1)$ & $(q^2{-}1)/2$ & -- \\ \hline \xrowht[()]{4pt}
$[a(\zeta^j),0],\ j\in {\cal V} $ & $(q^2-1)/8$ & $2q^2(q^2-1)$ & $(q^2{+}1)/2$ & -- \\ \hline
\hline
\end{tabular}
\caption{Table of conjugacy classes of untwisted elements of $G=M(q^2)$.}\label{tab:untwisted}
\end{table}

We proceed by determining conjugacy classes of twisted elements in $G=M(q^2)$. By \cite[Propositions 5 and 6]{EHS}, every twisted element $[A,1]\in G$ is conjugate in $G$ to an element of the form $[B,1]$ such that $B=\dia(\theta,1)$ or $B=\off(\theta,1)$ for some $\theta\in N(F)$. Conjugacy of twisted elements of $G$ in the overgroup $\ovl{G}$ was further investigated in detail in Proposition 7, 8 and Theorem 1 of \cite{EHS}. An inspection of the proofs of the three results with emphasis on comparison of conjugacy in $\ovl{G}$ and $G$ leads to the following.

\bp\label{prop:conj}
Let $\xi$ be a primitive element of $F$ and let $[A,1]$ be a twisted element of $G$. Then, exactly one of the following two cases occur:
\begin{itemize}
\item[{\rm 1.}] There exists exactly one odd $j\in \{1,2,\ldots,q-2\}$ such that $[A,1]$ is conjugate in $G$ to $[B,1]$ with $B=\dia(\xi^j,1)$, of order $2(q-1)/\gcd\{q-1,j\}$. The stabiliser of $[B,1]$ in $G$ is isomorphic to the cyclic group $C_{2(q-1)}$ generated by (conjugation by) a twisted element $[P,1]\in G$ with $P=\dia(\xi^{j\,'},1)$ for $j\,'=j-\frac{1}{2}(j-1)(q+1)$.
\item[{\rm 2.}] There exists exactly one odd $j\in \{1,2,\ldots,q\}$ such that $[A,1]$ is conjugate in $G$ to $[C,1]$ with $C=\off(\xi^j,1)$, of order $2(q+1)/\gcd\{q+1,j\}$. The stabiliser of $[C,1]$ in $G$ is isomorphic to the cyclic group $C_{2(q+1)}$ generated by (conjugation by) a twisted element $[Q,1]\in G$ with $Q=\off(\xi^{j\,'},1)$ for $j\,'=j+\frac{1}{2}(j-1)(q-1)$. \hfill
\end{itemize}
\ep

\pr As indicated, a proof can be obtained in an almost verbatim way from the statements and the proofs of Propositions 7 and 8, and Theorem 1, of \cite{EHS}. For readers interested in checking the details we mark here the differences in these proofs that are significant for distinction between conjugacy in $G$ and $\ovl{G}$.
\smallskip

Let $\xi$ be a primitive element of $F$. In the cases 1 and 2 of Proposition 7, twisted elements $[B,1]$ and $[B',1]$ with $B=\dia(\theta,1)$ and $B'=\dia(\theta',1)$ for $\theta,\theta'\in N(F)$ are conjugate in $\ovl{G}$ but {\sl not} in $G$ if, and only if, the ratio $\theta'/\theta$ in the case 1, and the product $\theta'\theta$ in the case 2, have the form $\xi^{t(q-1)}$ for $t$ odd. Similarly, in cases 1 and 2 of Proposition 8, elements $[B,1]$ and $[B',1]$ for $B=\off(\theta,1)$ and $B'=\off(\theta',1)$ are conjugate in $\ovl{G}$ but {\sl not} in $G$ if, and only if, $\theta'/\theta$ in the case 1, and $\theta'\theta$ in the case 2, have the form $\xi^{t(q+1)}$ for $t$ odd (i.e., they are non-squares in the unique subfield of $F$ of order $q$). All the remaining facts in the proofs of Propositions 7 and 8 apply to conjugacy in $G$.
\smallskip

Oddness of the values of $t$ above has the following two consequences in the proof of Theorem 1 of \cite{EHS} for restriction to conjugacy in $G$. Firstly, the congruences appearing in parts 1 and 2 of the proof have to be taken mod $2(q-1)$ and $2(q+1)$ and {\sl not} just mod $(q-1)$ and $(q+1)$, respectively, leading to upper bounds for $j$ in parts 1 and 2 of our Proposition that are two times larger than the bounds in \cite[Theorem 1]{EHS}. Secondly, parts 3 and 4 of the statement of Theorem 1 in \cite{EHS} refer to special cases arising due to the presence of exceptional conjugating elements in the corresponding parts of the proof. These need not be considered in our Proposition, because the exceptional conjugating elements turn out to lie outside $G$. Again, the remaining arguments apply to conjugacy in $G$. \hfill $\Box$
\medskip

By Proposition \ref{prop:conj}, there are a total of $q$ conjugacy classes of twisted elements in $G$, consisting of $(q-1)/2$ classes of `diagonal type' described in part 1, and $(q+1)/2$ classes of `off-diagonal type' from part 2, with centralizers of order $2(q-1)$ and $2(q+1)$, respectively. The explicit form of the twisted representatives from Proposition \ref{prop:conj} then immediately gives our Table \ref{tab:twisted} displaying conjugacy classes of twisted elements in $G$.
\smallskip

\begin{table}[ht]
	\centering
\begin{tabular}{|c||c|c|c|c|}
\hline \xrowht[()]{4pt}
Twisted representatives & \# Classes & $|{\rm Class}|$ & $|C_G|$ & Notes\\ \hline \hline \xrowht[()]{4pt}
$[\dia(\xi^j,1),1],\, {\rm odd}\ j\le q{-}2$ & $(q-1)/2$ & $q^2(q^2{+}1)(q{+}1)/2$ & $2(q{-}1)$ & -- \\ \hline \xrowht[()]{4pt}
$[\off(\xi^j,1),1],\, {\rm odd}\ j\le q$ & $(q+1)/2$ & $q^2(q^2{+}1)(q{-}1)/2$ & $2(q{+}1)$ & -- \\ \hline\hline
\end{tabular}
\caption{Table of conjugacy classes of twisted elements of $G=M(q^2)$.}\label{tab:twisted}
\end{table}

\section{Preliminary results on characters of $M(q^2)$}\label{sec:pre-char}

Summing up the results of Tables \ref{tab:untwisted} and \ref{tab:twisted}, the group $G=M(q^2)$ splits into a total of $(q+1)(q+5)/4$ conjugacy classes. Note that in $G$ the number of conjugacy classes of twisted elements, $q$, is by an order of magnitude smaller than the number of conjugacy classes of untwisted elements, $\approx q^2/4$, so that one may attempt to determine the character table of $G$ by lifting the character from its index-two subgroup $H\cong \PSL(2,q^2)$. The character table of $\PSL(2,q^2)$ is known and we will reproduce here a modification of the one from \cite[pp. 147--148]{JW}, see Table \ref{tab:PSL} below.

\begin{table}[ht]
	\centering
\begin{tabular}{|c||c|c|c|c|c|c|}
\hline \xrowht[()]{4pt}
$\PSL(2,q^2)$ & $I$ & $B_1$ & $B_{\veps}$ & $w$ & $a(\xi^j)$ & $a(\zeta^k)$  \\ \hline\hline \xrowht[()]{4pt}
$\#$ classes & $1$ & $1$ & $1$ & $1$ & $(q^2{-}5)/4$ & $(q^2{-}1)/4$ \\ \hline \xrowht[()]{4pt}
$|{\rm class}|$ & $1$ & $(q^4{-}1)/2$ & $(q^4{-}1)/2$ & $q^2(q^2{+}1)/2$ & $q^2(q^2{+}1)$ & $q^2(q^2{-}1)$ \\ \hline \xrowht[()]{4pt}
$|C_H|$ & $q^2(q^4{-}1)/2$ & $q^2$ & $q^2$ & $q^2{-}1$ & $(q^2{-}1)/2$ & $(q^2{+}1)/2$ \\ \hline\hline \xrowht[()]{4pt}
$\iota$ & $1$ & $1$ & $1$ & $1$ & $1$ & $1$ \\ \hline \xrowht[()]{4pt}
St & $q^2$ & $0$ & $0$ & $1$ & $1$ & $-1$ \\ \hline \xrowht[()]{4pt}
$\rho$ & $(q^2{+}1)/2$ & $(1{+}q)/2$ & $(1{-}q)/2$ & $1$ & $(-1)^j$ & $0$ \\ \hline \xrowht[()]{4pt}
$\rho'$ & $(q^2{+}1)/2$ & $(1{-}q)/2$ & $(1{+}q)/2$ & $1$ & $(-1)^j$ & $0$ \\ \hline \xrowht[()]{4pt}
$\rho_{\ell}$ & $q^2{+}1$ & $1$ & $1$ & $2(-1)^{\ell}$ & $\alpha^{j\ell} + \alpha^{-j\ell}$ & $0$ \\ \hline \xrowht[()]{4pt}
$\pi_m$ & $q^2{-}1$ & $-1$ & $-1$ & $0$ & $0$ & ${-}\beta^{km}{-}\beta^{-km}$ \\ \hline \hline
\end{tabular}
\caption{Table of irreducible characters of $H\cong \PSL(2,q^2)$.}\label{tab:PSL}
\end{table}

The symbols $B_1$, $B_{\veps}$, $w$, $a(\xi^j)$ and $a(\zeta^k)$ used in Table \ref{tab:PSL} for representatives of conjugacy classes are the same as explained at the beginning of Section \ref{sec:conj}, and $|C_H|$ is the order of the centralizer of the corresponding element in $H$. The irreducible characters are $\iota$ (the trivial one), St (the Steinberg permutation character), $\rho$, $\rho'$, $\rho_{\ell}$ for $1\le \ell \le (q^2-5)/4$, and $\pi_m$ for $1\le m\le (q^2-1)/4$. As before, $\xi$ and $\zeta$ are a primitive element of $F$ and a primitive $(q^2{+}1)^{\rm st}$ root of $1$ in $F'$, and the powers $j$ and $k$ in $\xi^j$ and $\zeta^k$ are bounded by $1\le j\le (q^2-5)/4$ and $1\le k\le (q^2-1)/4$. Finally, in a somewhat non-standard notation, $\alpha = \exp(4\pi i/(q^2{-}1))$ and $\beta = \exp(4\pi i/(q^2{+}1))$ are {\em complex} primitive roots of unity, respectively,of order $(q^2-1)/2$ and $(q^2+1)/2$.
\smallskip

Observe that $\PSL(2,q^2)$ has a total of $(q^2+5)/2$ conjugacy classes. This is roughly twice the number of conjugacy classes of $M(q^2)$ for large $q$. An explanation offered by the previous two sections is that $G$-conjugation fuses `most' pairs of $H$-conjugacy classes of untwisted elements in $H$ but there are only $q$ conjugacy classes of twisted elements of $G$.
\smallskip

The degrees of irreducible characters of $\PSL(2,q^2)$ follow from Table \ref{tab:PSL}, and those of $G\cong M(q^2)$ have been determined in \cite{Whi}. For convenience we display both in tabular form, with the proviso that degree $20$ cannot occur in the exceptional case when $q=3$:

\begin{table}[ht]
	\centering
\begin{tabular}{|c|c|}
\hline \xrowht[()]{4pt}
\ Degrees of irreducible characters of $H$ \ & \ Degrees of irreducible characters of $G$ \ \\ \hline \xrowht[()]{4pt}
$1$, \ $q^2$, \ $q^2-1$, \ $(q^2+1)/2$, \ $q^2+1$ & $1$, \ $q^2$, \ $2(q^2-1)$, \ $q^2+1$, \ $2(q^2+1)$ \\
\hline
\end{tabular}
\caption{Degrees of irreducible characters of $H\cong \PSL(2,q^2)$ and $G\cong M(q^2)$.}\label{tab:deg}
\end{table}

For reference to standard concepts and results in the theory of group characters we will use the monograph \cite{JL}. We will focus on results on characters of a group $G$ with a normal subgroup $H$ of index $2$; up to the end of the proof of Lemma \ref{lem:induc} the pair $G,H$ may be arbitrary but later we will return to our situation of $G$ and $H$ standing for $M(q^2)$ and $\PSL(2,q^2)$.
\smallskip

The {\em restriction} $\chi_H$ of a character $\chi$ of $G$ to $H$ is the character of $H$ defined by $\chi_H(g) =\chi(g)$ if $g\in H$ and $\chi_H(g)=0$ for $g\in G{\setminus} H$. The following is a short summary of results of \cite[Ch. 20]{JL} we need here; they only assume that $H$ is a normal subgroup of $G$ of index $2$. As usual, we let $\lambda$ denote the alternating character of $G$, with values $1$ on elements of $H$ and $-1$ on elements in $G{\setminus}H$.

\bl\label{lem:JaLi}
Let $\chi$ be an irreducible character of $G$. Then, $\chi_H$ is irreducible if and only if $\chi(g)\ne 0$ for some $g\in G{\setminus}H$, which is equivalent to $\chi\ne \chi\lambda$; moreover, in this case $\chi_H$ determines the pair $\{\chi, \chi\lambda\}$ uniquely. If $\chi_H$ is reducible, then it is the sum of two {\em distinct} irreducible characters of $H$ of the same degree. \hfill $\Box$
\el

As usual, the symbols $C_H(h)$ and $C_G(h)$ will denote the centralizers of an element $h\in H$ in $H$ and $G$. If $\varphi$ is a character of of our index-two subgroup $H$ of $G$, the corresponding {\em induced character} $\varphi^G$ of $G$ is given as follows (for some fixed $g\in G{\setminus}H$):
\begin{center}
$ \varphi^G(h) = \begin{cases}
\varphi(h) + \varphi(ghg^{-1}) & \mbox{if } h\in H \mbox{ and } C_H(h) = C_G(h); \\
2\varphi(h) & \mbox{if } h\in H \mbox{ and } C_H(h) \ne C_G(h); \\
0 & \mbox{if } h\in G{\setminus}H.
\end{cases}$
\end{center}
We note that if $C_H(h)=C_G(h)$, then conjugacy in $G$ fuses a pair of distinct $H$-conjugacy classes $C\ni h$ and $C'\ni ghg^{-1}$ in $H$ to a single $G$-conjugacy class (still in $H$); each such unordered pair $\{C,C'\}$ will be called a {\em fusion pair}.
\smallskip

For brevity we will refer to the value of the standard inner product $\langle \chi,\chi \rangle_G$ for a character $\chi$ of $G$ as the {\em norm} of $\chi$, denoted by $||\chi||_G$. A similar notation will be used for the norm of characters of $H$, and we will drop the subscript if the group is clear from the context. It is well known that $\chi$ is irreducible if and only if it has norm $1$. We will need the following auxiliary result on values of the norm of an induced character, and we only prove it for real characters (those with all values real), although the argument can easily be adapted to complex characters in general.

\bl\label{lem:induc}
Let $H$ be a normal subgroup of $G$ of index $2$ and let $\varphi$ be a real irreducible character of $H$. Then, $||\varphi^G||\in \{1,2\}$, with $||\varphi^G||=1$ if and only if $\varphi^G$ is an irreducible character of $G$, and $||\varphi^G||=2$ if and only if for every fusion pair $(C,C')$ of $H$-conjugacy classes in $H$ the values of $\varphi$ on $C\cup C'$ are constant.
\el

\pr
Let ${\cal F}$ and ${\cal F}'$ be the set of $H$-conjugacy classes in $H$ that belong, respectively, to some fusion pair and to no fusion pair of $H$. Since $\varphi$ is assumed to be real, the norm of $\varphi^G$ can be expressed in the form
\[ ||\varphi^G||_G = \frac{1}{|G|}\left( \sum_{C\in {\cal F}} \sum_{h\in C} (\varphi(h)+\varphi(ghg^{-1}))^2 + \sum_{C\in {\cal F}'} \sum_{h\in C} (2\varphi(h))^2  \right) \]
for any fixed $g\in G{\setminus}H$. The set ${\cal F}$ can be partitioned into fusion pairs $\{C,C'\}$; let ${\cal F}_0$ denote a subset of ${\cal F}$ consisting of $|{\cal F}|/2$ classes no two of which form a fusion pair. We will now use the fact that, for a fusion pair $\{C,C'\}$, the value of $\varphi(h)+\varphi(ghg^{-1})$ for $h\in C$ is the same as the value of $\varphi(h')+\varphi(gh'g^{-1})$ for $h'=g^{-1}hg\in C'$ (note that $g^2\in H$). Using this (and $|G|=2|H|$) the above expression for $||\varphi^G||$ can be rewritten as follows, with the first sum being taken over conjugacy classes in ${\cal F}_0$:
\be\label{eq:norm} ||\varphi^G||_G = \frac{1}{|H|}\left( \sum_{C\in {\cal F}_0} \sum_{h\in C} (\varphi(h)+\varphi(ghg^{-1}))^2 + 2\sum_{C\in {\cal F}'} \sum_{h\in C} \varphi(h)^2  \right) \ee
With the help of the obvious inequality $(x+y)^2 \le 2(x^2+y^2)$ for real $x,y$ (with equality if and only if $x=y$) we obtain from \eqref{eq:norm} the inequality
\be\label{eq:norm2} ||\varphi^G||_G \le \frac{1}{|H|}\left( \sum_{C\in {\cal F}_0} \sum_{h\in C} 2(\varphi(h)^2 +\varphi(ghg^{-1})^2) + 2\sum_{C\in {\cal F}'} \sum_{h\in C} \varphi(h)^2  \right) = 2||\varphi||_H \ee
with equality if and only if $\varphi(h)=\varphi(ghg^{-1})$ for every $h$ that belongs to a conjugacy class forming a fusion pair.
\smallskip

Since norms are positive integers and the character $\varphi$ of $H$ was assumed to be irreducible, that is, $||\varphi||_H=1$, from \eqref{eq:norm2} we obtain $||\varphi^G||_G\in \{1,2\}$. Moreover, the above necessary and sufficient condition for equality in \eqref{eq:norm2}, that is, for $||\varphi^G||_G = 2$, translates into the condition that for every fusion pair $(C,C')$ of $H$-conjugacy classes the values of $\varphi$ on $C$ are the same as the values of $\varphi$ on $C'$. \hfill $\Box$
\medskip

The following observation will also be useful; from this point on we will return to our notation $G=M(q^2)$ with a subgroup $H\cong \PSL(2,q^2)$ of index $2$.

\bl\label{lem:useful}
Let $\varphi$ be an irreducible character of $H\cong \PSL(2,q^2)$ such that $\varphi^G = \chi + \chi'$ for two irreducible characters of $G=M(q^2)$. Then, $\chi'=\chi\lambda \ne \chi$.
\el

\pr Lemma \ref{lem:JaLi} implies that the induced characters $\chi^{}_H$ and $\chi'_H$ of $H$ are either both reducible, or both irreducible. We begin by eliminating the first possibility.
\smallskip

Suppose that both $\chi^{}_H$ and $\chi'_H$ are reducible characters of $H$. By inspecting the values in Table \ref{tab:deg} we find that that their degrees {\em must} be $q^2+1$. Since they are restrictions of irreducible characters, it follows that $\chi$ and $\chi'$, respectively, are obtained from $\chi^{}_H$ and $\chi'_H$ by setting their values to be zero everywhere in $G{\setminus}H$. By Lemma \ref{lem:JaLi} we have $\chi^{}_H=\psi + \eta$ and $\chi'_H=\psi'+\eta'$ for some collection of irreducible characters of $H$ of degree $(q^2+1)/2$ each, with $\psi\ne\eta$ and $\psi'\ne\eta'$. However, by Table \ref{tab:PSL} the group $H$ has only two irreducible characters of such a degree, so that $\chi^{}_H=\chi'_H$, and as the restricted characters uniquely determine $\chi$ and $\chi'$ we obtain $\chi=\chi'$. Now, $||\varphi^G||_G=2$ whereas $||2\chi||_G = 4$, a contradiction.
\smallskip

Thus, both $\chi^{}_H$ and $\chi'_H$ are irreducible characters of $H$. Since $\varphi^G$ is zero on $G{\setminus}H$, from $||\varphi^G||_G=2$ it follows that $||(\varphi^G)_H||_H=4$. But $||(\varphi^G)_H||_H = \langle \chi^{}_H+\chi'_H, \chi^{}_H+ \chi'_H\rangle_H$, which is only equal to $4$ if $\chi^{}_H=\chi'_H$. Lemma \ref{lem:JaLi} then implies that $\chi'= \chi\lambda \ne\chi$. \hfill $\Box$
\medskip

\section{Lifting characters of $\PSL(2,q^2)$ to $M(q^2)$}\label{sec:lift-char}

We begin this section by identifying the irreducible characters of $G$ of degree $1$ and $q^2$, the unique corresponding characters of $H$ being the trivial character $\iota$ and the Steinberg character ${\rm St}$. Recall that the (irreducible) Steinberg character of a finite $2$-transitive permutation group is evaluated at any permutation in the group by subtracting $1$ from the number of fixed points of the permutation. The standard actions of the groups $H\cong \PSL(2,q^2)$ and $G=M(q^2)$, on the $q^2+1$ projective points are both $2$-transitive, and even sharply $3$-transitive in the case of $G$. It is easy to verify that the values of the Steinberg character of $G$, which we will denote {\scalebox{0.6}{$\Sigma$}}, are the same as those of ${\rm St}$ on conjugacy classes of untwisted elements, and $+1$ and $-1$ on conjugacy classes of twisted diagonal and twisted off-diagonal elements, respectively.

\bl\label{lem:lift-2}
Let $\varphi$ be an irreducible character of $H$. If $\deg(\varphi)=1$, then $\varphi^G = \iota + \lambda$, and $\iota$ with $\lambda$ are the only irreducible characters of $G$ of degree $1$. If $\deg(\varphi)=q^2$, then one has $\varphi^G =$ {\scalebox{0.6}{$\Sigma$}} $+$ {\scalebox{0.6}{$\Sigma$}}\,$\lambda$, where {\scalebox{0.6}{$\Sigma$}} is the Steinberg permutation character, and {\scalebox{0.6}{$\Sigma$}} with {\scalebox{0.6}{$\Sigma$}}\,$\lambda$ are the only irreducible characters of $G$ of degree $q^2$.
\el

\pr By Table \ref{tab:deg} in both cases we have $\varphi^G=\chi+\chi'$ for irreducible characters $\chi$ and $\chi'$ of $G$, so that Lemma \ref{lem:useful} applies and $\varphi^G=\chi+\chi\lambda$. The rest is a consequence of uniqueness of irreducible characters of $H$ of degree $1$ and $q^2$.  \hfill $\Box$
\medskip

With the help of Lemma \ref{lem:lift-2} and known facts from character theory we are in position to determine the {\em number} of irreducible characters of $G$ of degrees appearing Table \ref{tab:deg}.

\bl\label{lem:deg}
The number of irreducible characters of $G=M(q^2)$ of a given degree are as follows:
\el

\vskip 1mm
\
\vskip -11mm

{\rm
\begin{table}[ht]
	\centering
\begin{tabular}{|c||c|c|c|c|c|}
\hline \xrowht[()]{4pt}
Character degree & $1$ & $q^2$ & $2(q^2-1)$ & $q^2+1$ & $2(q^2+1)$ \\ \hline \xrowht[()]{4pt}
$\#$ irreducible characters & $2$ & $2$ & $(q^2-1)/8$ & $2q-3$ & $(q-1)(q-3)/8$ \\ \hline
\end{tabular}
\end{table}
}
\medskip

\pr By Lemma \ref{lem:lift-2} we know that the number of irreducible characters of $G$ of degree $1$ and $q^2$ is $2$ in both cases. Further, by Table \ref{tab:deg} there are only three remaining degrees of irreducible characters of $G$, namely, $2(q^2-1)$, $q^2+1$ and $2(q^2+1)$; let $c_1,c_2,c_3$, respectively, be the numbers of such characters. Now, $2+2+c_1+c_2+c_3 = (q+1)(q+5)/4$ is the number of conjugacy classes in $G$, and as the sum of squares of all character degrees is equal to $|G|^2$, we also have
\[ 2\times 1^2 + 2\times q^4 + c_1\times 4(q^2-1)^2 + c_2\times (q^2+1)^2 + c_3\times 4(q^2+1)^2 = q^2(q^4-1) \]
which can be shown to be equivalent to
\be\label{eq:sq}
((c_2+4c_3)(q^2+1) + 2)(q^2+1) = (q^2-1)^2(q^2-4c_1)
\ee
Since $(q^2-1)/2$ and $(q^2+1)/2$ are relatively prime and the second one is odd, by the factorization appearing in \eqref{eq:sq} the number $(q^2+1)/2$ must divide the (odd) number $q^2-4c_1$, and so $q^2+1$ is a divisor of $2(q^2-4c_1)$. This, however, is possible only if the two numbers are equal, that is, $q^2+1 = 2(q^2-4c_1)$, which is if and only if $c_1=(q^2-1)/8$. Having determined the value of $c_1$ we are left with a system of two equations in two unknowns and it is easy to check that its unique solution if $c_2=2q-3$ and $c_3=(q-1)(q-3)/8$. \hfill $\Box$
\medskip

We are in position to give substantial information about lifts of irreducible characters of $H$ onto $G$. Its statement includes remarks on integrality of values of the characters, and refers to the notation used in the table of conjugacy classes of untwisted elements of $M(q^2)$ (Table \ref{tab:untwisted}, in particular the sets ${\cal U}$ and ${\cal V}$) and in Table \ref{tab:PSL} of characters of $\PSL(2,q^2)$; the first two assertions from Lemma \ref{lem:lift-2} are included for completeness.

\bp\label{prop:lift}
Let $\varphi$ be an irreducible character of $H\cong \PSL(2,q^2)$ with the induced character $\varphi^G$ of $G \cong M(q^2)$. Then, exactly one of the following cases occur:
\begin{itemize}
\item[{\rm (1)}] $\deg(\varphi^G)=2$ and $\varphi^G = \iota + \lambda$, with $\iota$ and $\lambda$ being the only irreducible characters of $G$ of degree $1$; they have integral values and are non-zero on $G{\setminus}H$;
\item[{\rm (2)}] $\deg(\varphi^G)=2q^2$ and $\varphi^G =$ {\scalebox{0.6}{$\Sigma$}} $+$ {\scalebox{0.6}{$\Sigma$}}\,$\lambda$, where {\scalebox{0.6}{$\Sigma$}} is the Steinberg permutation character, and {\scalebox{0.6}{$\Sigma$}} with {\scalebox{0.6}{$\Sigma$}}\,$\lambda$ are the only irreducible characters of $G$ of degree $q^2$; they are again integral and have non-zero values on $G{\setminus}H$;
\item[{\rm (3)}] $\deg(\varphi^G)=2(q^2-1)$ and $\varphi^G$ is one of the $(q^2-1)/8$ irreducible characters of $G$ of degree $2(q^2-1)$ such that $(\varphi^G)_H = \pi_m + \pi_{mq}$ for $m\in {\cal V}$; these characters $\varphi^G$ are all real and with all-zero values on $G{\setminus}H$;
\item[{\rm (4)}] $\deg(\varphi^G)=q^2+1$ and $\varphi^G$ is a unique irreducible character of $G$ of degree $q^2+1$ with $(\varphi^G)_H=\rho+\rho'$; it is integral and identically zero on $G{\setminus}H$;
\item[{\rm (5)}] $\deg(\varphi^G)=2(q^2+1)$ and $\varphi^G$ is one of the $(q-1)(q-3)/8$ irreducible characters of $G$ of degree $2(q^2+1)$ such that $(\varphi^G)_H = \rho_\ell + \rho_{\ell q}$ for $\ell\in {\cal U}$; these $\varphi^G$ are all real and all-zero on $G{\setminus}H$;
\item[{\rm (6)}] $\deg(\varphi^G)=2(q^2+1)$ and $\varphi^G = \chi + \chi\lambda$ for $q-2$ pairs $\chi\ne\chi\lambda$ of irreducible characters of $G$ of degree $q^2+1$, with $\chi_H= \rho_\ell$ for the $(q-3)/2$ values $\ell=r(q+1)/2$ such that $1\le r\le (q-3)/2$, and the $(q-1)/2$ values  $\ell=s(q-1)/2$, $1\le s\le (q-1)/2$; there is a total of $2q-4$ such distinct irreducible characters of $G$ and each of these have at least one non-zero value on $G{\setminus}H$.
\end{itemize}
\ep

\pr
As noted, we may skip the first two items, and among the remaining ones we begin with (4). If $\deg(\varphi^G)= q^2+1$, then $\varphi^G$ is necessarily an irreducible character of $G$ such that $(\varphi^G)_H$ is the sum of two distinct irreducible characters of $H$ of degree $(q^2+1)/2$. But by Table \ref{tab:PSL} there are only two such characters of $H$, namely, $\rho$ and $\rho'$, and one may check that they both induce the {\em same} character of $G$, giving the conclusion of (4).
\smallskip

Next, we turn our attention to (6), where $\varphi^G$ is assumed to be reducible and of degree $2(q^2+1)$. By Lemma \ref{lem:useful} we have $\varphi^G=\chi + \chi\lambda$ for some irreducible character $\chi$ of $G$ of degree $q^2+1$, with $\chi_H$ also irreducible. This means that $\chi_H$ must be one of the characters $\rho_\ell$ for a suitable $\ell$, and the same applies to $\varphi$, of course. But $||\varphi^G|| =2$, so that by Lemma \ref{lem:induc} the character $\varphi=\rho_\ell$ must be constant on any fusion pair of $H$-conjugacy classes. By the findings of section \ref{sec:conj} applied to this case, the $H$-conjugacy classes generated by the elements $[a(\xi^j),0]$ and $[a(\xi^{jq}),0]$ for $j\in \{1,\ldots,(q^2-5)/4\}$ form a fusion pair if they are distinct. The previous condition therefore means that the values of $\rho_\ell$ and $\rho_{\ell q}$ on every such pair of classes must be the same (so that we may ignore distinctness here). By Table \ref{tab:PSL} this translates, for a fixed $\ell$, to the equality of the complex numbers $\alpha^{j\ell} + \alpha^{-j\ell}$ and $\alpha^{j\ell q} + \alpha^{-j\ell q}$ for {\em all} $j$ as above.
\smallskip

A simple calculation (as in section \ref{sec:conj}) reveals that the two complex numbers coincide if and only if $(\alpha^{j\ell(q+1)}-1) (\alpha^{j\ell(q-1)}-1)=0$ (or, equivalently, $\alpha^{j\ell q}\in \{\alpha^{j\ell},\alpha^{-j\ell}\}$). To fulfil the condition on the constant value on fusion classes we are looking for the values of $\ell\in \{1,\ldots,(q^2-5)/4\}$ for which the last equation holds for {\em every} $j\in \{1,\ldots, (q^2-5)/4\}$, which happens if and only if $\ell$ is a multiple of $(q+1)/2$ or $(q-1)/2$. In our range $1\le\ell\le (q^2-5)/4$ given by Table \ref{tab:PSL} this yields the $q-2$ values of $\ell$ appearing in the statement of (6), and hence also the $q-2$ possibilities for $(\varphi^G)_H=\rho_\ell + \rho_{\ell q} = 2\rho_\ell$. The latter give $2(q-2)$ irreducible characters $\chi$ and $\chi\lambda$ of degree $q^2+1$, all with some non-zero value on $G{\setminus}H$, and such that $\chi_H=(\chi\lambda)_H=\rho_\ell$. Since the total number of irreducible characters of $G$ of degree $q^2+1$ is $2q-3$ by Lemma \ref{lem:deg} and one of such characters was identified in part (4), the total number of irreducible characters of $G$ referred to in (6) is $2q-4$, as claimed.
\smallskip

We continue with (5), assuming that $\varphi^G$ is irreducible and has degree $2(q^2+1)$. The restriction of $\varphi$ to $H$ cannot be irreducible, so that $(\varphi^G)_H$ is a sum of two distinct irreducible characters of $H$ of degree $q^2+1$. Such characters are all of the form $\rho_\ell$ for suitable values of $\ell$, and $\varphi$ in this case must also be equal to one of these. Further, our assumption on the degree of $\varphi^G$ together with Lemma \ref{lem:induc} imply that there must be at least one fusion pair of $H$-conjugacy classes such that $\varphi$ is not constant on the pair. The fusion pairs here are generated by the same classes as in (6), and our calculations in the proof of (6) imply that a fusion pair on which $\varphi$ is not constant exists if and only if $\ell$ is {\em not} one of the values listed in the conclusion of (6). One may check that this reduces to the condition $\ell\in {\cal U}$ for the set ${\cal U}$ introduced in section \ref{sec:conj}, with $(\varphi^G)_H=\rho_\ell + \rho_{\ell q}$ for $\ell\in {\cal U}$. (We note again that all the accompanying calculations are analogous to those in section \ref{sec:conj} where the condition $\theta^q \notin \{\pm \theta, \pm\theta^{-1}\}$ was considered. Here the condition on $\alpha$ translates to $\alpha^{j\ell q}\notin \{\alpha^{j\ell},\alpha^{-j\ell}\}$ and leads to the same conclusion that $\ell\in {\cal U}$ because the expression $(\varphi^G)_H=\rho_\ell + \rho_{\ell q}$ is invariant under the substitutions $\ell\mapsto -\ell$ and $\ell\mapsto \ell q$; the $\pm 1$ term is absorbed by the factor $4$ in $\alpha=\exp(4\pi i/(q^2-1))$.)
\smallskip

By Table \ref{tab:PSL} it may be verified that if $\varphi = \rho_\ell$ for $\ell\in {\cal U}$, then $(\varphi^G)_H = \rho_\ell + \rho_{\ell q}$, so that $\rho_\ell$ and $\rho_{\ell q}$ determine the {\em same} lift onto $G$. Referring again to Tables \ref{tab:untwisted} and \ref{tab:PSL} one obtains this way a total of $|{\cal U}|= (q-1)(q-3)/8$ irreducible characters $\varphi^G$ of $G$ whose restriction to $H$ has the form $\rho_\ell + \rho_{\ell q}$ for $\ell\in U$. All such lifts are distinct, and by Lemma \ref{lem:deg} there cannot be any other irreducible character of $G$ of degree $2(q^2+1)$.
\smallskip

Finally, let us consider (3), assuming that $\deg(\varphi^G)=2(q^2-1)$. By Table \ref{tab:deg}, $\varphi^G$ must be an irreducible character of $G$, and as its restriction to $H$ cannot be irreducible, $(\varphi^G)_H$ is a sum of two distinct irreducible characters of $H$ of degree $q^2-1$. The latter are all of the form $\pi_m$ for suitable values of $m$, and $\varphi$ must also be of this form. In section \ref{sec:conj} we showed that the $H$-conjugacy classes of the elements $[a(\zeta^j),0]$ and $[a(\zeta^{jq}),0]$ form a fusion pair for every $j\in {\cal V}$. One may check that the values of any given $\pi_m$ are non-constant on at least one of these fusion pairs, which conforms to the last part of Lemma \ref{lem:induc}. Further, a direct verification against Table \ref{tab:PSL} shows that if $\varphi = \pi_m$, then $(\varphi^G)_H = \pi_m + \pi_{mq}$. Thus, both $\pi_m$ and $\pi_{mq}$ determine the {\em same} lift onto $G$, and by Tables \ref{tab:untwisted} and \ref{tab:PSL} they give rise to $(q^2-1)/8$ irreducible characters of $G$ of the form $\pi_m + \pi_{mq}$ corresponding to the values $m\in {\cal V}$. It may be verified that these lifts are distinct, and by Lemma \ref{lem:deg} there is no other irreducible character of $G$ of degree $2(q^2-1)$.
\smallskip

As regards remarks on integral and real values, most of them are obvious, and note that $\varphi^G$ is always identically zero on $G{\setminus}H$. This completes the proof. \hfill $\Box$
\medskip

An inspection of the numbers in Proposition \ref{prop:lift} reveals that in the cases (1) -- (5) there are, respectively, $2+2+(q^2-1)/8+1+(q-1)(q-3)/8$ {\em real} characters, and all of them have been completely determined. This makes a total of $(q^2-2q+21)/4$ irreducible characters of $G$ that are all real. Subtracting this from the number of all irreducible characters of $G$ leaves us only with the $2q-4$ ones described in part (6) of Proposition \ref{prop:lift}, that may assume complex values, but only on $G{\setminus}H$ as all the characters of $H$ are real.
\smallskip

We now address the question of possible character values that are not real. A well known result in the theory of group characters is that, for an element $g$ of a group, the values of all irreducible characters of the group evaluated at $g$ are real if and only if $g$ is conjugate to $g^{-1}$ in the group, or, equivalently, the conjugacy class of $g$ is closed under inversion in the group; such conjugacy classes are called {\em real}. By Table \ref{tab:untwisted} every conjugacy class of untwisted elements in $G=M(q^2)$ is real. For conjugacy of twisted elements we have, in the notation of Proposition \ref{prop:conj}:

\bl\label{lem:real}
The conjugacy class of the twisted element $[\dia(\xi^j,1),1]$ in $G=M(q^2)$ for odd $j$, $1\le j\le q-2$, is real if and only if $q\equiv 3$ {\rm mod $4$} and $j=(q-1)/2$. The conjugacy class of $[\off(\xi^j,1),1]$ for odd $j$, $1\le j\le q$, is real if and only if $q\equiv 1$ {\rm mod $4$} and $j=(q+1)/2$.
\el

\pr Let $d(j)=[\dia(\xi^j,1),1]$ and $d'(j)=[\off(\xi^j,1),1]$. Observe first that $d(j)^{-1}=d(-jq)$ and $d'(j)^{-1}=d'(jq)$. By Proposition \ref{prop:conj} the elements $d(j)$ and $d(-jq)$ are in the same conjugacy class in $G$ for odd $j$, $1\le j\le q-2$, if and only if $j\equiv -jq$ mod $2(q-1)$. This can only happen if $q\equiv 3$ mod $4$, and then the congruence is equivalent to $j(q+1)/4 \equiv 0$ mod $(q-1)/2$. As the modulus of the last congruence is relatively prime to $(q+1)/4$ it follows that this only leaves us with the value $j=(q-1)/2$ in our range for $j$. The analysis for the elements $d'(j)$ and $d'(jq)$ is analogous.  \hfill $\Box$
\medskip

This gives, in both cases mod 4, exactly $q-1$ conjugacy classes of twisted elements in $G$ that are not real. By the earlier remarks we also know that there are at most $2q-4$ irreducible characters of $G$ (those from part (6) of Proposition \ref{prop:lift}) that are not real.

\section{Representations for the remaining characters}

In the previous section we have almost completely determined the character table of $G=M(q^2)$ and we have been left with $2q-4$ `missing' irreducible characters, which are the only ones that may assume non-real values on $q-1$ conjugacy classes of $G$ formed by elements of $G{\setminus}H$. These characters have been referred to in part (6) of Proposition \ref{prop:lift} and from now on we will denote them by $\chi_\ell$ and $\chi_\ell \lambda$, with $(\chi_\ell)_H =(\chi_\ell \lambda)_H=\rho_\ell$, for the total of $q-2$ values of $\ell$ in the set $L= L^+\cup L^-$, where $L^+= \{r(q+1)/2;\ 1\le r\le (q-3)/2\}$ and $L^-= \{s(q-1)/2;\ 1\le s\le (q-1)/2\}$. We will determine these characters in the next section; here we first derive some related representations by lifting one-dimensional representations of a suitable subgroup of $H$. The method is an adaptation of derivation of principal series representations for two-dimensional special linear groups (cf. \cite[p. 232]{Dor}).
\smallskip

Let $H_{\rm upp}$ be the subgroup of $H<G$, $H\cong \PSL(2,q^2)$, induced by upper-triangular matrices with determinant $1$. Explicitly, if $\xi\in F^*=F{\setminus}\{0\}$ is a fixed primitive element as before and if $h(u,d)$ is the $2\times 2$ matrix with first and second row of the form $(\xi^u,d)$ and $(0,\xi^{-u})$ for an arbitrary non-negative integer $u < (q^2-1)/2$ and any $d\in F$, then
\be\label{eq:upper}
H_{\rm upp} = \{ [h(u,d),0]\in G;\ 0\le u < (q^2{-}1)/2,\ d\in F \}\ ;
\ee
observe that $|H|=q^2(q^2-1)/2$. It is well known (see e.g. \cite[p. 232]{Dor} adapted to the projective case) that for every $\ell\in L$ (and, in fact, for every integer $\ell$ but this will not be needed here) the assignment
\be\label{eq:phiC} \Phi_\ell:\ [h(u,d),0] \mapsto \exp\left(\frac{4\pi i}{q^2-1}\ell u\right)
\ee
defines a one-dimensional complex representation of $H_{\rm upp}$. To lift such a representation onto one of the entire group $G= M(q^2)$ we proceed as ibid by first constructing a suitable set of coset representatives of $H_{\rm upp}$ in $G$. To do this, for $0\le t < q^2-1$ we introduce $2\times 2$ lower-triangular matrices $m(t)$ and $m'(t)$, as well as matrices $m(\infty)$ and $m'(\infty)$, as follows:
\be\label{eq:mm'}
m(t) {=} \left( \begin{array}{cc} 1 & 0 \\ \xi^t & 1\end{array}\right), \
m'(t) {=} \left( \begin{array}{cc} \xi & 0 \\ \xi^t & 1\end{array}\right), \
m(\infty) {=} \left( \begin{array}{cc} 0 & 1 \\ -1 & 0\end{array}\right), \
m'(\infty) {=} \left( \begin{array}{cc} 0 & \xi \\ -1 & 0\end{array}\right).
\ee
Using the $2q^2$ matrices from \eqref{eq:mm'} we introduce $2(q^2+1)$ elements of $G$ by letting
\begin{eqnarray*}
&x_t=[m(t),0], \ {\rm and}\ y_t=[m'(t),1] \ {\rm for} \ 0\le t < q^2-1,\ {\rm together\ with} \\ &x_\infty = [m(\infty),0], \ \ y_\infty=[m'(\infty),1], \ \ x_\ast=[I,0],\ {\rm and}\  y_\ast=[\dia(\xi,1),1];
\end{eqnarray*}
note that $x_\ast$ is the unit element of $G$. It may be checked that this set of $n_q=2(q^2+1)$ elements of $G$ is a left transversal for the subgroup $H_{\rm upp}$.
\smallskip

With the help of this transversal we will now lift any one-dimensional representations $\Phi\in \{\Phi_\ell;\ \ell\in L\}$ of $H_{\rm upp}$ described in \eqref{eq:phiC} onto an $n_q$-dimensional representation $\Phi^G$ of $G$; the method of lifting or inducing originates from \cite{Fr2}. Before doing so we will make an agreement about indexation. Let ${\rm Ind}=\{0,1,\ldots, q^2-2,\ast, \infty\}$, where the entries $0,1,\ldots,q^2-2$ are considered mod $q^2-1$, and let ${\rm Ind}'=\{z';\ z\in {\rm Ind}\}$. To describe $n_q\times n_q$ matrices we will use the $n_q$ indices from the set ${\rm Ind}\cup {\rm Ind}'$ equipped with the linear ordering
\be\label{eq:order}
0< 1< \ldots< (q^2-2)< 0'< 1'< \ldots< (q^2-2)'< \ast< \infty< \ast'< \infty'\ .
\ee

Invoking now \cite[Lemma 9.1]{Dor} adapted to our situation and using the introduced notation, an $n_q$-dimensional representation $\Phi^G$ is obtained by assigning, to every element $g\in G$ the $n_q\times n_q$ matrix $\Phi^G(g)$ whose $(a,b)$-th entry for $a,b\in {\rm Ind}\cup {\rm Ind}'$ is determined by the following rules (using the convention that $(z')'=z$ for our indices):
\begin{center}
$ \Phi^G(g)_{a,b} = \begin{cases}
\ \Phi(x_agx_b^{-1}) & \mbox{if } g\in H,\ a,b\in {\rm Ind} \ \mbox{ and } x_agx_b^{-1}\in H_{\rm upp}; \\
\ \Phi(y_{a'}gy_{b'}^{-1}) & \mbox{if } g\in H,\ a,b\in {\rm Ind}' \ \mbox{ and } y_{a'}gy_{b'}^{-1}\in H_{\rm upp}; \\
\ \Phi(x_agy_{b'}^{-1}) & \mbox{if } g\in G{\setminus}H,\   a\in {\rm Ind},\ b\in {\rm Ind}' \ \mbox{ and } x_agy_{b'}^{-1}\in H_{\rm upp}; \\
\ \Phi(y_{a'}gx_b^{-1}) & \mbox{if } g\in G{\setminus}H,\  a\in {\rm Ind}',\ b\in {\rm Ind} \ \mbox{ and } y_{a'}gx_b^{-1}\in H_{\rm upp}; \\
\ 0 & \mbox{ in all other cases }.
\end{cases}$
\end{center}

By part (6) of Proposition \ref{prop:lift}, each lifted representation $\Phi^G$ for $\Phi\in \{\Phi_\ell;\ \ell\in L^+\cup L^-\}$ is reducible and splits into two irreducible representations of dimension $q^2+1$ each. The missing irreducible characters $\chi_\ell$ and $\chi_\ell\lambda$ on $G{\setminus}H$ for $\ell\in L^+\cup L^-$ are given by traces of these representations, or, equivalently, by traces corresponding to the two $G$-invariant $(q^2+1)$-dimensional subspaces of the representation $\Phi^G$.
\smallskip

For evaluation of these traces (which we will do in the next section) we will need explicit knowledge of $\Phi^G$-images of a generating set of $G=M(q^2)$ and its  subgroup $H\cong \PSL(2,q^2)$. It is well known that $H\cong \PSL(2,q^2)$ is generated by the two elements $g_1=[\dia(\xi,\xi^{-1}),0]$ and $g_2=[\dia(1,0)+ \off(-1,1),0]$, of order $(q^2-1)/2$ and $3$, respectively. In accordance with Table \ref{tab:twisted}, as representatives of conjugacy classes of $M(q^2){\setminus}H$ we will take the elements $g_3=g_3(j)= [\dia(\xi^j,1),1]$ and $g_4=g_4(j) = [\off(\xi^j,1),1]$ for odd positive integers $j\le q{-}2$ and $j\le q$, respectively. For the $\Phi^G$-images of these elements we obtain:

\bp\label{prop:miss-rep}
For every $\ell\in L$ the lifts $\Phi_\ell^G$ of the linear representation $\Phi_\ell$ evaluated at the elements $g_1$, $g_2$, $g_3$ and $g_4$ are $n_q\times n_q$ unitary matrices with entries as follows, where $\gamma_\ell = \exp(2\pi i \ell/(q^2-1))=\exp(\pi i r/(q-1))$, and indices in ${\rm Ind}\cup {\rm Ind}'$ distinct from $\ast,\ast',\infty,\infty'$ are understood mod $q^2-1$:
\medskip

\noindent {\rm (a)} If $g=g_1=[\dia(\xi,\xi^{-1}),0]$, then
\begin{center}
$ \Phi_\ell^G(g)_{a,b} = \begin{cases}
\ \gamma_\ell^2 & \mbox{if } (a,b)\in \{(\ast,\ast),\ (t,t{+}2),\ 0\le t\le q^2-2\}, \\
\ \gamma_\ell^{2q} & \mbox{if } (a,b)\in \{(\ast',\ast'),\ (t',(t{+}2q)'),\ 0\le t\le q^2-2\}, \\
\ \gamma_\ell^{-2} & \mbox{if }  (a,b)=(\infty,\infty),  \\
\ \gamma_\ell^{-2q} & \mbox{if }  (a,b)=(\infty',\infty'),  \\
\ 0 & \mbox{in all other cases }.
\end{cases}$
\end{center}

\noindent {\rm (b)} If $g=g_2=[\dia(1,0)+\off(-1,1),0]$, then, letting $\ovl{q}=(q^2-1)/2$,
\begin{center}
$ \Phi_\ell^G(g)_{a,b} = \begin{cases}
\ \gamma_\ell^{-2a} & \mbox{if }  (a,b)=(t,f(t)) \ {or}\ (t',f(t)'), 0\le t\le q^2-2\ {and}\ t\ne \ovl{q}, \\
\ 1 & \mbox{if } (a,b)\ {or}\ (a',b')\ {is\ in\ the\ set}\ \{(\ovl{q},\ast),(\ast,\infty),(\infty,\ovl{q})\}, \\
\ 0 & \mbox{in all the remaining cases, }
\end{cases}$
\end{center}
where the function $f$ on residue classes $t$ mod $q^2-1$ and $t\ne \ovl{q}$ is given by $\xi^{-t}+\xi^{f(t)}=-1$.\\

\noindent {\rm (c)} If $g=g_3=g_3(j)=[\dia(\xi^j,1),1]$ for an odd positive integer $j\le q-2$, then
\begin{center}
$ \Phi_\ell^G(g)_{a,b} = \begin{cases}
\ \gamma_\ell^{j-1} & \mbox{if } a=t \ {\rm and}\ b=(t{+}j)', \ 0\le t\le q^2-2, \ {\rm or} \  (a,b)=(\ast,\ast'), \\
\ \gamma_\ell^{jq+1} & \mbox{if } a=t' \ {\rm and}\ b=t{+}jq, \ 0\le t\le q^2-2, \ {\rm or} \  (a,b)=(\ast',\ast), \\
\ \gamma_\ell^{-j-1} & \mbox{if } (a,b)=(\infty,\infty'),\\
\ \gamma_\ell^{-jq+1} & \mbox{if } (a,b)=(\infty',\infty),\\
\ 0 & \mbox{in all other cases }.
\end{cases}$
\end{center}

\noindent {\rm (d)} If $g=g_4=g_4(j)=[\off(\xi^j,1),1]$ for an odd positive integer $j\le q$, then
\begin{center}
$ \Phi_\ell^G(g)_{a,b} = \begin{cases}
\ (-1)^\ell\gamma_\ell^{-j-2t-1} & \mbox{if } a=t \ {\rm and}\ b=({-}t{-}j)', \ 0\le t\le q^2-2, \\
\ (-1)^\ell\gamma_\ell^{-jq-2t+1} & \mbox{if } a=t' \ {\rm and}\ b={-}t{-}jq, \ 0\le t\le q^2-2, \\
\ (-1)^\ell\gamma_\ell^{j-1} & \mbox{if } (a,b)=(\ast,\infty'), \\
\ (-1)^\ell\gamma_\ell^{-j-1} & \mbox{if } (a,b)=(\infty,\ast'), \\
\ (-1)^\ell\gamma_\ell^{jq+1} & \mbox{if } (a,b)=(\ast',\infty), \\
\ (-1)^\ell\gamma_\ell^{-jq+1} & \mbox{if } (a,b)=(\infty',\ast),\\
\ 0 & \mbox{in all the remaining cases }.
\end{cases}$
\end{center}

\ep

\pr
For illustration we will only deal with the entire part (b) and the second items of both (c) and (d), as the verification in all the remaining cases is analogous; the matrices are unitary by inspection and hence so is the entire lifted representation $\Phi_\ell^G$. For all calculations we note that $\Phi_\ell^G$ was introduced before the statement of Proposition \ref{prop:miss-rep}, preceded by an exposition of the associated matrices in \eqref{eq:mm'} together with the elements $x_a$ and $y_b$.
\smallskip

{\sl Part (b).}  Letting $g=g_2=[\dia(1,0)+\off(-1,1),0]$, for $a,b\in \{0,1,\ldots,t^2-2\}$ one has $\Phi_\ell^G(g)_{a,b}=\Phi(x_agy_{b'}^{-1})$ if $x_agy_{b'}^{-1} \in H_{\rm upp}$; otherwise $\Phi_\ell^G(g)_{a,b}=0$. By the rules (explained in Section \ref{sec:pre}) for calculations in our group $G$, the inverse of $y_t=[m'(t),1]$ is $y_t^{-1}=[n(t)^\sigma,1]$ where $n(t)$ is the $2\times 2$ matrix with rows $(\xi^{-1},0)$ and $(-\xi^{t-1} ,1)$. Following these rules (and slightly abusing the notation and using $a,b\notin \{\ast,\infty\}$ also as exponents at $\xi$), the product $X=x_agx_{b}^{-1}$ evaluates to
\[ X=\left(\begin{array}{cc} 1 & 0 \\ \xi^a & 1\end{array}\right)
\left(\begin{array}{cc} 1 & -1 \\ 1 & 0\end{array}\right)
\left(\begin{array}{cc} 1 & 0 \\ -\xi^{-b} & 1\end{array}\right) =
\left(\begin{array}{cc} \xi^b+1 & -1 \\ \xi^a+\xi^{a+b}+1 & -\xi^a\end{array}\right)\ .\]
The condition $X\in H_{\rm upp}$ is equivalent to $\xi^{-a}+\xi^b=-1$, which, for $a\ne \ovl{q}$, determines $b$ as a function $f(a)$ from the equation $\xi^{-a}+\xi^{f(a)}=-1$. In such a case the main diagonal of $X$ consists of the elements $1+\xi^{f(a)}=-\xi^{-a}$ and $-\xi^a$, so that by \eqref{eq:phiC} the value of $\Phi_\ell^G(g)_{a,f(a)}$ is equal to $\gamma_\ell^{-2a}$. One may check that for the pair $(a',b')$, $0\le a,b\le q^2-2$, evaluation of the product $y_agy_{b}^{-1}$ reduces to determining membership of $m'(a)gn(b)$ in  $H_{\rm upp}$ and gives the same condition on $a,b$ as above, that is, $b=f(a)$ for $a\notin \{\ovl{q},\ast,\infty\}$, with the same value $\gamma_\ell^{-2a}$ of $\Phi_\ell^G(g)_{a',f(a)'}$. Calculating all of $x_{\ovl{q}}g x_\ast^{-1}$, $x_\ast gx_\infty^{-1}$, $x_\infty gx_\ast^{-1}$ and their $y$-versions one obtains matrices with main diagonal $1,1$, implying the second item of (b).
\smallskip

Regarding the function $f$, mutual equivalence of the three equations $\xi^{-a}+\xi^b+1=0$, $\xi^{-b}+\xi^{-a-b}+1=0$ and $\xi^{a+b}+\xi^a+1=0$ implies that $f(a)=b$, $f(b)=-a-b$ and $f(-a-b)=a$ for $a\notin  \{\ovl{q},\ast,\infty\}$. Since raising in $q$-th power is an automorphism of $F={\rm GF}(q^2)$ one also has $f(qa)=qf(a)$. It follows that if $q$ is not a power of $3$ then $f$ as a permutation of the undashed indices not in $\{\ovl{q},\ast,\infty\}$ has two fixed points, namely, the two primitive $3$rd roots of $1$ in $F$, and for the remaining values, $f$ consists of $3$-cycles of the form $(a,b,c)$ where $b=f(a)$ and $a+b+c=0$; if $q$ is a power of $3$ then $f$ has $0$ as its unique fixed point and the remaining orbits of $f$ are $3$-cycles as above. An analogous statement applies to dashed indices.
\smallskip

{\sl The second item of part (c).} For $g=g_3(j)=[\dia(\xi^j,1),1]$ finding the value of $y_{a'}gx_b^{-1}$ by the calculation rules in $G$ reduces to evaluating the product $X=m'(a)\cdot\dia(\xi^{jq},1) \cdot m(b)^{-1}$, which results in
\[ X=\left(\begin{array}{cc} \xi & 0 \\ \xi^a & 1\end{array}\right)
\left(\begin{array}{cc} \xi^{jq} & 0\\ 0 & 1\end{array}\right)
\left(\begin{array}{cc} 1 & 0 \\ -\xi^b & 1\end{array}\right) =
\left(\begin{array}{cc} \xi^{jq+1} & 0 \\ \xi^{jq+a}-\xi^b & 1 \end{array}\right) \ .\]
Here one has $X\in H_{\rm upp}$ if and only if $b=a+jq$ mod $(q^2-1)$. To evaluate $\Phi_\ell$ at such an $X$ one needs to represent it in the form $[h(u,d),0]$, which reduces to looking for an element $z\in F$ such that $zX$ has a pair of mutually inverse entries $\xi^u$ and $\xi^{-u}$ in the main diagonal. The obvious choice here is $z= \xi^{-(jq+1)/2}$, giving by \eqref{eq:phiC} the value $\Phi_\ell^G(g)_{a',a+jq} = \gamma_\ell^{jq+1}$.
\smallskip

{\sl The second item of part (d).} Here, in $x_agy_{b'}^{-1}= [m(a)\cdot\off(\xi^j,1) \cdot n(b),0]$ the product $X$ of the three matrices is
\[ X=\left(\begin{array}{cc} 1 & 0 \\ \xi^a & 1\end{array}\right)
\left(\begin{array}{cc} 0 & \xi^j \\ 1 & 0\end{array}\right)
\left(\begin{array}{cc} \xi^{-1} & 0 \\ -\xi^{b-1} & 1\end{array}\right) =
\left(\begin{array}{cc} -\xi^{b+j-1} & \xi^j \\ \xi^{-1}-\xi^{a+b+j-1} & \xi^{a+j}\end{array}\right).\]
Now, $X$ belongs to $H_{\rm upp}$ if and only if $b=-a-j$ mod $(q^2-1)$, and to determine the corresponding value of $\Phi_\ell$ we again need to represent this element as $[h(u,d),0]$ by finding some $z\in F$ such that the main diagonal of $zX$ has the form $(\xi^u,\xi^{-u})$. Using $b+j=-a$ and $(-1)= \xi^{(q^2-1)/2}$, multiplication by $z=\xi^{-(2j-2 +q^2-1)/4}$ produces the required diagonal entries in $zX$ for $u=-a-(j+1)/2 - (q^2-1)/4$. By \eqref{eq:phiC}, the value of $\Phi_\ell$ at $[h(u,d),0]\in H_{\rm upp}$ for $h(u,d)=zX$ is then equal to $\exp(4\pi i \ell u/(q^2-1))$. Finally, with the help of $\gamma_\ell = \exp(2\pi i \ell/(q^2-1))$ and the substitution for $u$ the value of $\Phi_\ell$ simplifies to $(-1)^\ell\gamma_\ell^{-j-2a-1}$, which is as claimed if $a=t$ and $b=({-}t{-}j)'$.
\hfill $\Box$
\medskip

It will be of advantage to restate the above results in form of block matrices. For $\ell\in L$ let $P=P_\ell$ be the matrix of dimension $q^2+1$ indexed by the ordered set ${\rm Ind}$ with $P_{\infty,\infty} = \gamma_\ell^{-1}$ and $P_{\ast,\ast} = P_{a,a{+}1} = \gamma_\ell$ for every $a\in {\rm Ind}{\setminus} \{\ast,\infty\}$ mod $(q^2-1)$, and with all the remaining entries equal to zero. We will leave out the subscript $\ell$ if no confusion is likely. Further, let $A_{\rm upp}$, $A_{\rm low}$, $B$, $C_{\rm upp}=C_{\rm upp}(j)$ and $C_{\rm low}= C_{\rm low}(j)$ for odd positive $j$ be square matrices of dimension $q^2+1$, indexed by the ordered set Ind, and defined as follows:
\begin{eqnarray}\label{eq:ABC}
& A_{\rm upp}=P^2, \ \ \ A_{\rm low}=P^{2q}= A_{\rm upp}^q  \nonumber \\
& B_{a,f(a)}=\gamma_\ell^{-2a} {\ \rm if\ } a\in {\rm Ind}{\setminus} \{\ovl{q},\ast,\infty\},\ B_{a,b}=1 {\ \rm if\ } (a,b)\in\{(\ovl{q},\ast),(\ast,\infty), (\infty,\ovl{q})\} \ \ \ \ \ \\
& C_{\rm upp} = C_{\rm upp}(j)=\gamma_\ell^{-1}P^j,\ \ C_{\rm low}= C_{\rm low}(j)=\gamma_\ell P^{jq}
\nonumber
\end{eqnarray}
where entries of $B$ not listed are assumed to be zero; the parameter $j$ will be omitted if no loss of clarity is likely. From now on we will also extend the symbols $\dia$ and $\off$ also to $2\times 2$ blocks of dimension $q^2+1$. One may check that, in this new notation, the matrices $A^{\dagger}=\Phi_\ell^G(g_1)$, $B^{\dagger}=\Phi_\ell^G(g_2)$ and $C^\dagger = C^{\dagger}(j)=\Phi_\ell^G (g_3(j))$ from Proposition \ref{prop:miss-rep} can be displayed in the form
\be\label{eq:blockform}
A^{\dagger} = \dia(A_{\rm upp},A_{\rm low}),  \ \ B^{\dagger} = \dia(B,B), \ \ C^{\dagger} = \off(C_{\rm upp}, C_{\rm low})\
\ee
where we assume a natural extension of indexation from Ind to Ind' in the bottom $q^2+1$ coordinates of the `larger' matrices of dimension $2(q^2+1)$; this also justifies the usage of the subscripts `upp' and `low' for the upper and lower non-zero blocks of the matrices of \eqref{eq:ABC} which appear in \eqref{eq:blockform}.
\smallskip

It is well known that the group $H\cong \PSL(2,q^2)$ is generated by the pair of elements $g_1$ and $g_2$. Since the matrices $A^{\dagger}$ and $B^{\dagger}$ of dimension $2(q^2+1)$ assigned to $g_1$ and $g_2$ in the representation $\Phi_\ell^G$ are block-diagonal, it follows that {\sl all} matrices of $\Phi_\ell^G$ assigned to elements of $H$ are block-diagonal (with diagonal blocks of dimension $q^2+1$). Thus, the restriction $\Phi_\ell^H$ of the unitary representation $\Phi_\ell^G$ to $H$ is a sum of two irreducible $(q^2+1)$-dimensional unitary representations, say, $\Phi^{(1)}$ and $\Phi^{(2)}$, of the group $H$; formally, $\Phi_\ell^H=\Phi^{(1)}\oplus \Phi^{(2)}$.
\smallskip

This is a good point to recall that our aim is to determine the $2(q-2)$ characters denoted $\chi_\ell$ and $\chi_\ell \lambda$ at the beginning of the previous section (which are the characters of part (6) of Proposition \ref{prop:lift}) for $\ell\in L$ with $|L|=q-2$. As the restrictions $(\chi_\ell)_H$ and $(\chi_\ell \lambda)_H$ coincide and are equal to the character $\rho_\ell$ of $H\cong \PSL(2,q^2)$, it follows that in the sum $\Phi_\ell^H=\Phi^{(1)}\oplus \Phi^{(2)}$ the constituents $\Phi^{(1)}$ and $\Phi^{(2)}$, generated, respectively, by the pairs $A_{\rm upp},B$ and $A_{\rm low},B$, must be {\sl equivalent} unitary irreducible representations of $H$. This implies the existence of an `intertwining' unitary matrix $M=M_\ell$ of dimension $q^2+1$ with the property that $\Phi^{(1)}(h)M = M\Phi^{(2)}(h)$ for every $h\in H$. In particular, using the common notation $Y^M=M^{-1}YM$ for conjugation, for the constituents of $A^\dagger$ and $B^\dagger$ in \eqref{eq:blockform} one has
\be\label{eq:AMBM}
A_{\rm low} = A_{\rm upp}^q = A_{\rm upp}^M \ \ \ {\rm and}  \ \ \ B=B^M \ .
\ee
Further, using $A_{\rm upp}=P^2$, $A_{\rm low}=P^{2q} = A_{\rm upp}^q$ from \eqref{eq:ABC} and realizing that $P^{q^2}=P$ one obtains
\[ A_{\rm upp}^{M^2} = A_{\rm low}^M = (A_{\rm upp}^q)^M = (A_{\rm upp}^M)^q =  A_{\rm low}^q = A_{\rm upp}^{q^2} = A_{\rm upp} \ .\]
Since from \eqref{eq:AMBM} we obviously have $B^{M^2}=B$, it follows that $M^2$ commutes with the representation $\Phi^{(1)}$. By Schur's theorem, $M^2$ is a constant multiple of the identity matrix. However, for the purpose of intertwining the constant multiple may be an arbitrary non-zero complex number, which we will henceforth choose to be equal to $1$. Thus, without loss of generality we may assume that $M^2=I$, which means that $M$ is a unitary {\sl Hermitian} matrix; now we also have ${\rm det}(M)=\pm 1$.
\smallskip

\section{Determination of the remaining characters}

We begin this section by determining all the non-trivial $G$-invariant (and hence $(q^2+1)$-dimensional) subspaces of our $2(q^2+1)$-dimensional representation $\Phi_\ell^G$. The method relies on the following result which extends \cite[Lemma 2.6, p. 56]{TB} and is likely to be folklore in representation theory; we therefore include only a short proof.

\bl\label{lem:invar}
Let $\Psi_1\oplus\Psi_2$ be a direct sum of a pair of equivalent complex irreducible representations of a group $K$ with associated disjoint vector spaces $V_1$ and $V_2$ and with an intertwining matrix $N$. If a non-trivial subspace $V$ of $V_1\oplus V_2$ is $K$-invariant, then either the projection of $V$ onto exactly one of $V_1$, $V_2$ is zero, or there is a non-zero $c\in \mathbb{C}$ such that $V=\{({\bf u},{\bf u}(cN));\ {\bf u}\in V_1\}$.
\el

\pr By default, all of $\Psi_i$ and $V_i$ for $i=1,2$ as well as $V$ and $N$ must have the same dimension. Letting $\pi_i$ be the projection of $V_1\oplus V_2$ onto $V_i$ for $i=1,2$, by irreducibility the intersection $V\cap {\rm ker}(\pi_i)$ is either $V_i$ or trivial. Leaving the possibility when (exactly) one of $\pi_i(V)$ is equal to $V_i$ we are left with the case when both intersections $V\cap {\rm ker}(\pi_i)$ are trivial. But then, letting $\pi_{i,V}$ denote the restriction of $\pi_i$ to $V$, it follows that $V_i=\pi_{i,V}(V)\cong V$ for $i=1,2$; in particular, the mappings $\pi_{i,V}$ are $K$-invariant isomorphisms. The composition $\Psi=\pi^{-1}_{1,V}\pi^{}_{2,V}$ is a $K$-invariant isomorphism $V_1\to V_2$, which, by Schur's theorem, must be given by ${\bf u}\Psi = {\bf u}(cN)$ for an intertwining matrix $N$, unique up to a non-zero multiplicative constant. Now, if ${\bf v}\in V$ is an arbitrary vector, letting ${\bf u} = {\bf v}\pi_{1|V}$ we obtain ${\bf v}=({\bf v}\pi_{1|V}, {\bf v}\pi_{2|V})= ({\bf u},{\bf u}\Psi) = ({\bf u},{\bf u}(cN))$.
\hfill $\Box$
\medskip

Applying Lemma \ref{lem:invar} to our restricted representation $\Phi_\ell^H$ of the subgroup $H<G$ one sees that an invariant subspace $W$ for $\Phi_\ell^H$ has either zero projection onto the first (undashed) $q^2+1$ coordinates in the set ${\rm Ind}$, or a zero projection onto the second (dashed) coordinates in ${\rm Ind}'$, or else has the form $W=\{({\bf u},c{\bf u}M);\ {\bf u}\in \mathbb{C}^{q^2+1}\}$, where $M = M_\ell$ is the intertwining involutory matrix of \eqref{eq:AMBM}. One of these then {\sl must} form an invariant subspace for the entire representation $\Phi_\ell^G$ of the group $G$. Realizing that the representation $\Phi_\ell^G$ is obtained from $\Phi_\ell^H$ by adjoining any of the matrices $C^\dagger$ from \eqref{eq:blockform} with zero diagonal blocks, the first two possibilities are immediately ruled out, and we obtain:

\bc\label{cor:invar} Every non-trivial invariant subspace of $\Phi_\ell^G$ has the form $W=\{({\bf u},c{\bf u}M);$ ${\bf u}\in \mathbb{C}^{q^2+1}\}$, where $M$ is the intertwining involutory matrix of the representations $\Phi^{(1)}$ and $\Phi^{(2)}$ of $H$ generated by $A_{\rm upp},B$ and $A_{\rm low},B$, and $c$ is a non-zero complex constant. \hfill $\Box$
\ec

By the Corollary, for every vector ${\bf w}\in W$, i.e., of the form ${\bf w}=({\bf u},c{\bf u}M)$ for ${\bf u}\in \mathbb{C}^{q^2+1}$, the vector ${\bf w}C^{\dagger} = ({\bf u},c{\bf u}M) \off(C_{\rm upp}, C_{\rm low}) = (c({\bf u}M)C_{\rm low}, {\bf u}C_{\rm upp})$ must also belong to $W$. This means that ${\bf u}C_{\rm upp} = {\bf u}(cM)C_{\rm low}(cM)$, and as ${\bf u}$ was arbitrary and $M=M^{-1}$, it follows that
\be\label{eq:MCM} C_{\rm low} = c^{-2}C_{\rm upp}^M \ .
\ee

With the help of Corollary \ref{cor:invar} and the knowledge of $C_{\rm upp} = \gamma_\ell^{-1}P^j$, $C_{\rm low}= \gamma_\ell P^{jq}$ and $A_{\rm upp}=P^2$ by \eqref{eq:ABC} we now determine the possible values of the constant $c$. Namely, squaring the expressions for $C_{\rm upp}$ and $C_{\rm low}$ one obtains
\[ C_{\rm upp}^2 = \gamma_\ell^{-2}P^{2j} = \gamma_\ell^{-2}A_{\rm upp}^j \ \ \ {\rm and} \ \ \ C_{\rm low}^2 = \gamma_\ell^{2}P^{2jq} = \gamma_\ell^{2}A_{\rm upp}^{jq}\ ,  \]
and substituting these into the square of \eqref{eq:MCM} gives
\[  \gamma_\ell^{2}A_{\rm upp}^{jq} = c^{-4} \gamma_\ell^{-2}(A_{\rm upp}^j)^M\ .\]
The last equation reduces by \eqref{eq:AMBM}, i.e., by $A_{\rm upp}^M=A_{\rm upp}^q$, to $c^4 = \gamma_\ell^{-4}$, so that $c= \pm i^{\delta(\ell)}\gamma_\ell^{-1}$ for some $\delta(\ell)\in \{0,1\}$, where $i$ is the complex imaginary unit. Hence, \eqref{eq:MCM} can equivalently be written in the form
\be\label{eq:MCM1} C_{\rm low} = (-1)^{\delta(\ell)}\gamma_\ell^2C_{\rm upp}^M \ .
\ee

Observe that if $c$ is one of the four values determined above for which the subspace $W=\{({\bf u},c{\bf u}M);\ {\bf u}\in \mathbb{C}^{q^2+1}\}$ is $G$-invariant, then $W^\perp=\{({\bf v},-c{\bf v}M);\ {\bf v}\in \mathbb{C}^{q^2+1}\}$ is another such subspace and the pair $(W,W^\perp)$ forms an orthogonal  decomposition of $\mathbb{C}^{2(q^2+1)}$. Indeed, letting ${}^*$ denote the complex conjugate transpose and using the fact that $M$ is unitary ($MM^*=I$) together with $cc^* =1$, evaluating the standard inner product of (complex) vectors from $W$ and $W^\perp$ gives
\[ ({\bf u},c({\bf u}M))\cdot ({\bf v},-c({\bf v}M)) = {\bf u}{\bf v}^* - c{\bf u}M((c{\bf v}M)^* = {\bf u}{\bf v}^* - cc^*{\bf u}MM^*{\bf v}^* =0 \ .\]

Consider such a pair $W,W^\perp$ of $G$-invariant subspaces. From the way our lifted representation $\Phi_\ell^G$ was introduced before the statement of Proposition \ref{prop:miss-rep} and from the calculations in its proof it follows that for every $g\in G{\setminus}H$ the unitary matrix $\Phi_\ell^G(g)$ has the form $E^\dagger=E^\dagger(g)= \off(E_{\rm upp},E_{\rm low})$ with off-diagonal (and necessarily unitary) blocks of dimension $q^2+1$ each. Consider a complex eigenvalue $x$ of $E^\dagger$ associated with a non-zero eigenvector ${\bf w}$, splitting uniquely as ${\bf w}= {\bf w}_1+{\bf w}_2$ for ${\bf w}_1\in W$ and ${\bf w}_2\in W^\perp$. Now, $({\bf w}_1+{\bf w}_2)E^\dagger={\bf w}E^\dagger=x{\bf w}=x({\bf w}_1+{\bf w}_2)$ and by $\mathbb {C}^{2(q^2+1)} = W\oplus W^\perp$ and the $G$-invariance of the two subspaces we obtain ${\bf w}_1E^\dagger=x{\bf w}_1$ and ${\bf w}_2E^\dagger=x{\bf w}_2$. It follows that the eigenspace ${\rm Eig}(E^\dagger,x)$ of $E^\dagger$ associated with the eigenvalue $x$ is a direct sum $({\rm Eig}(E^\dagger,x)\cap W) \oplus ({\rm Eig}(E^\dagger,x)\cap W^\perp)$ of the corresponding spaces of eigenvectors of $E^\dagger$ that belong to $W$ and $W^\perp$. By diagonalizability of unitary matrices this also means that the spectrum of $E^\dagger$, of size $2(q^2+1)$, is a concatenation of the spectra ${\rm Spec}(E^\dagger,W)$ and ${\rm Spec}(E^\dagger,W^\perp)$ of size $q^2+1$ consisting, respectively, of the eigenvalues corresponding to eigenspaces contained in $W$ and $W^\perp$.
\smallskip

Let now ${\bf w}=({\bf u},c{\bf u}M)\in W$ be an eigenvector of $E^\dagger$ for an eigenvalue $x$. Since $W$ is $G$-invariant and hence preserved by $E^\dagger$, one has
\[ (x{\bf u},xc{\bf u}M) = x{\bf w} = {\bf w}E^\dagger = (c{\bf u}ME_{\rm low},{\bf u}E_{\rm upp}) \in W \ .\]
Membership of the last vector in $W$ means that ${\bf u}E_{\rm upp} = c{\bf u}ME_{\rm low}(cM)$, and since this holds for {\sl every} eigenvector ${\bf u}$ corresponding to any eigenvalue $x$, we obtain $E_{\rm upp}=c^2E_{\rm low}^M$, or, equivalently, $c^*E_{\rm upp}M=cME_{\rm low}$. The chain of displayed equations further gives $ x{\bf u} = {\bf u}c ME_{\rm low}$ and $xc{\bf u}M = (c{\bf u}M)cE_{\rm low}M$, which means that ${\bf u}$ and $c{\bf u}M$ are eigenvectors of the matrices $cME_{\rm low} = c^*E_{\rm upp}M$ and $cE_{\rm low}M$, respectively, for the same eigenvalue $x$. Conversely, if ${\bf u}$ is an eigenvector of $c^*E_{\rm upp}M = cME_{\rm low}$ for an eigenvalue $x$, then right multiplication by $cM$ shows that $c{\bf u}M$ is an eigenvector of $cME_{\rm low}cM = E_{\rm upp}$ for $x$ and hence $x({\bf u},c{\bf u}M) = (c{\bf u}ME_{\rm low}, {\bf u}E_{\rm upp})$, that is, ${\bf w}=({\bf u},c{\bf u}M)\in W$ is an eigenvector of $E^\dagger$ for the same eigenvalue $x$.
\smallskip

This correspondence between the eigenspaces of $E^\dagger$ that are subspaces of $W$ and eigen\-spaces of the (unitary and hence diagonalizable) matrix $cME_{\rm low} = c^*E_{\rm upp}M$ together with the earlier established facts about spectra lead to the conclusion that the multi-set ${\rm Spec}(E^\dagger,W)$ is equal to the spectrum of $cME_{\rm low} = c^*E_{\rm upp}M$. This way we have arrived at the following conclusion for our missing character $\chi_\ell$ for $\ell\in L$, which we state as a summary of the above considerations together with \eqref{eq:AMBM} and \eqref{eq:MCM1}.

\bp\label{prop:eigM} Let $M$ be a unitary Hermitian matrix of dimension $q^2+1$ with determinant $\pm 1$ such that
\be\label{eq:ABCL+} A_{\rm low} = A_{\rm upp}^M\ , \ \ B=B^M \ \ {\rm and} \ \ C_{\rm low} = c^{-2}C_{\rm upp}^M = (-1)^{\delta(\ell)}\gamma_\ell^2C_{\rm upp}^M\ \ee
for some $c\in \{\pm i^{\delta(\ell)}\gamma_\ell^{-1}\}$. For an arbitrary $g\in G{\setminus}H$ let $\Phi_\ell^G(g) = \off(E_{\rm upp}(g),E_{\rm low}(g))$. Then, the pair of characters $\chi_\ell$ and $\chi_\ell\lambda$ for each $\ell\in L$ is determined by letting $\chi_\ell(g) = {\rm tr}(cME_{\rm low}(g)) ={\rm tr}(c^* E_{\rm upp}(g)M)$ for every $g\in G{\setminus}H$. \hfill $\Box$ \ep

Proposition \ref{prop:eigM} allows for a quick determination of the missing pair of characters $\chi_\ell$ and $\chi_\ell\lambda$ for $\ell \in L^+ = \{r(q+1)/2;\ 1\le r\le (q-3)/2\}$. The key observation now is that for $\ell= r(q+1)/2\in L^+$ ($1\le r\le (q-3)/2$) one has $\gamma_\ell^{q-1} = (-1)^r$ by Proposition \ref{prop:miss-rep}. Let $Q$ be the permutation matrix corresponding to the permutation of the set Ind fixing $\ast$ and $\infty$ and sending $a\mapsto aq$ for every $a\in {\rm Ind}{\setminus} \{\ast,\infty\}$ mod $(q^2-1)$. Observe that, due to $\gamma_\ell^{q-1} = (-1)^r$, for the matrix $P$ introduced immediately after the end of the proof of Proposition \ref{prop:miss-rep} one has $P^q=(-1)^rP^Q$. With this in hand and using oddness of $j$, an inspection of \eqref{eq:ABC} shows that
\be\label{eq:ABCL++} A_{\rm low} = A_{\rm upp}^Q\ , \ \ B=B^Q \ \ {\rm and} \ \ C_{\rm low} = (-1)^r\gamma_\ell (P^j)^Q = (-1)^r\gamma_\ell^2C_{\rm upp}^Q\ .\ee
Comparing \eqref{eq:ABCL++} with \eqref{eq:ABCL+} implies that for $\ell\in L^+$ one can simply take $M=Q$ for the intertwining matrix, with $(-1)^{\delta(\ell)}=(-1)^r$ and $c^{-1}=c^*=i^r\gamma_\ell$.
\smallskip

Applying Proposition \ref{prop:eigM} further, for $\ell\in L^+$ the value of the `missing' character $\chi_\ell$ at the element $g=g_3(j)$ of part (c) of Proposition \ref{prop:miss-rep} may be taken to be the trace of the product $c^*C_{\rm upp}(j)Q$. By \eqref{eq:ABC} one has $C_{\rm upp}(j) = \gamma_\ell^{-1}P^j$ and one may check that the only non-zero diagonal elements of the product $C_{\rm upp}(j)Q$ are those in positions $\ast$ and $\infty$, which are $\gamma_\ell^{j-1}$ and $\gamma_\ell^{-j-1}$. Therefore $\chi_\ell (g_3(j))= c^*(\gamma_\ell^{j-1} +\gamma_\ell^{-j-1}) = i^r(\gamma_\ell^{j} +\gamma_\ell^{-j})$ for $\ell= r(q+1)/2\in L^+$, $1\le r\le (q-3)/2$. By the same token, letting $\Phi_\ell^G(g_4(j)) = \off(D_{\rm upp},D_{\rm low})$ for the element $g_4(j)$ of part (d) of Proposition \ref{prop:miss-rep}, for $\ell\in L^+$ one has $\chi_\ell(g_4(j)) = {\rm tr}(c^*D_{\rm upp}(j)Q)$ and as this matrix has a zero diagonal by inspection it follows that $\chi_\ell(g_4(j)) = 0$.
\smallskip

Using the notation $\tau^+_\ell(j) = i^r(\alpha_r^j+\alpha_r^{-j})$ for $\ell= r(q+1)/2$, with $\alpha_r= \exp(i\pi r/(q{-}1)) = \gamma_\ell$, the missing pair of characters $\chi_\ell$, $\chi_\ell\lambda$ may be given as follows.

\bc\label{cor:L+} For $\ell= r(q+1)/2\in L^+$, $1\le r\le (q-3)/2$, the pair of characters $\chi_\ell$ and $\chi_\ell\lambda$ are determined by  $\chi_\ell(g_3(j)) = \tau^+_\ell(j)$ for odd positive $j\le q-2$, and $\chi_\ell(g_4(j)) = 0$ for odd positive $j\le q$. \hfill $\Box$
\ec

We continue by calculating the values of $\chi_\ell$ for $\ell\in L^-$ on conjugacy classes of twisted elements of $G$. Out of the previous results it is easy to extract, for each of the $(q+1)/2$ values of $\ell\in L^-$, the $(q-1)/2$-dimensional vector of values of $(\chi_\ell(g_3(j);\ {\rm odd}\ j\le q-2)$. Observe that, for any fixed $\ell\in L^-$, the restriction $(\chi_\ell)_H$ is orthogonal to {\em each} of the $(q-3)/2$ restrictions $(\chi_\ell)_H$ for $\ell\in L^+$ and to the trivial character $\iota$. It follows that, for our fixed $\ell\in L^-$, the $(q-1)/2$-dimensional vector ${\bf w}=(\chi_\ell(g_3(j);\ {\rm odd}\ j\le q-2)$ is orthogonal to the system of $(q-1)/2$ mutually orthogonal vectors consisting of the all-one vector and the $(q-3)/2$ vectors $(\gamma_\ell^j + \gamma_\ell^{-j};\ {\rm odd}\ j\le q-2)$ for the $(q-3)/2$ values of $\ell\in L^+$, all of dimension $(q-1)/2$. But such a ${\bf w}$ must then be the zero vector, which implies that $\chi_\ell$ is zero at all elements $g_3(j)$ for every $\ell\in L^-$.
\smallskip

It remains to determine the values of $\chi_\ell$ for $\ell\in L^-$ at the elements $g_4(j)$ for odd positive $j\le q$, which will take considerably more space. To this end it is sufficient to determine the matrix $M$ for this situation, and although general methods for calculating intertwining matrices are available (see e.g. \cite{MSH}) we adopt here a more direct approach. Summing up some of the facts established so far, the unitary Hermitian matrix $M$ with ${\rm det}(M)= \pm 1$ is, up to the sign, determined by the equations from Proposition \ref{prop:eigM}, the second and third of which are equivalent to $BM=MB$ and $(-1)^{\delta(\ell)}\gamma_\ell^2C_{\rm upp}M = MC_{\rm low}$. Written in terms of coordinates $M_{a,b}$ for $a,b\in {\rm Ind}$, the first equation gives
\be\label{eq:MB}
 M_{f(a),f(b)} = \gamma_\ell^{2(|a|-|b|)}M_{a,b}\ \ {\rm with} \ \ |a| = \begin{cases}
\ a \ \ {\rm if} \ a\ne \ast,\infty \\
\ 0 \ \ otherwise. \end{cases}
\ee
To deal with the second equation we will assume that $\ell=s(q-1)/2\in L^-$, $1\le s\le (q-1)/2$, and we will choose $\delta(\ell)$ such that $(-1)^{\delta(\ell)} =(-1)^s$; as we shall see, this choice will be consistent with the forthcoming calculations. Taking this into the account and realizing that now $\gamma_\ell^{q+1}=(-1)^s$, it can be checked that the equation $(-1)^{\delta(\ell)}\gamma_\ell^2C_{\rm upp}M = MC_{\rm low}$, taken first for $j=1$ and then extended by induction for every $t$ mod $q^2+1$, translates into the following linear system, where $a,b\notin \{\ast,\infty\}$:
\begin{align}\label{eq:MC}
& &M_{a+t,b+tq} &= \gamma_\ell^{-2t} M_{a,b} \ &\ \ \ \ \ \  \nonumber \\
\ \ \ \ \ M_{a+t,\ast} &= \gamma_\ell^{-2t} M_{a,\ast}, &M_{\ast,b+tq} &= \gamma_\ell^{-2t} M_{\ast,b}, &M_{\ast,\ast}=0 \ &\ \ \ \ \ \ \\
\ \ \ \ \ M_{a+t,\infty} &= M_{a,\infty}, &M_{\infty,b+tq} &= M_{\infty,b}, &M_{\infty,\infty}=0  \ &\ \ \ \ \ \ \nonumber
\end{align}
It can also be verified that the system \eqref{eq:MC} for $t=2$ implies the first equation $A_{\rm low} = A_{\rm upp}^M$ from Proposition \ref{prop:eigM}.
\smallskip

In what follows we will use the fact that $\ovl{q}=(q^2-1)/2$ is the only non-trivial involution mod $(q^2-1)$, with $q\ovl{q}\equiv \ovl{q}$ mod $(q^2-1)$; we will also frequently use the cycle $(\ovl{q},\ast,\infty)$ of the permutation $f$. Let $y=M_{0,\infty}$. The leftmost equation in the last row of \eqref{eq:MC} shows that $M_{a,\infty}=y$ for every $a\in {\rm Ind}_0={\rm Ind}{\setminus} \{\ast,\infty\}$. Applying first \eqref{eq:MB} with $\gamma_\ell^{2q}= \gamma_\ell^{-2}$, followed by using the middle equation of \eqref{eq:MC} for $t=q$ and finally by \eqref{eq:MB} again one obtains
\[ M_{\infty,f(\ovl{q}+1)}=\gamma_\ell^{-2}M_{\ast,\ovl{q}+1} = \gamma_\ell^{-2}(\gamma_\ell^2 M_{\ast,\ovl{q}}) = M_{\ovl{q},\infty} = y\ .\] But from $M=M^*$ it also follows that $M_{\infty,f(\ovl{q}+1)}=(M_{f(\ovl{q}+1),\infty})^* = y^*$, which in combination with the above implies $y=y^*$, that is, $y$ is a real number. Using \eqref{eq:MB} it also follows that $y=M_{\ovl{q},\infty}=M_{\infty,\ast}=M_{\ast,\infty}^*=M_{\ast,\infty}$, so that the entries of the column of $M$ marked $\infty$ are constantly equal to $y$ except for $M_{\infty,\infty}=0$. Due to $M=M^*$ the same conclusion is valid for the row of $M$ marked $\infty$. By $MM^*=M^2=I$, the dot product of the row and the column of $M$ marked $\infty$ must be equal to $1$, that is, $q^2y^2=1$, from which $y=\pm q^{-1}$.
\smallskip

The entries of $M$ in the row and column marked $\ast$ are determined by the equations of \eqref{eq:MC} containing asterisks. To determine the remaining entries, observe first that, by \eqref{eq:MC} and \eqref{eq:MB}, one has $M_{\infty,\infty}=M_{\ovl{q},\ovl{q}}=M_{0,0}=0$, and hence by the first equation of \eqref{eq:MC} one also has $M_{t,tq}=M_{tq,t}=0$ for every $t$ mod $(q^2-1)$. The remaining entries of $M$ turn out to be non-zero and can be determined as follows. For every $a'\in {\rm Ind}{\setminus} \{\ovl{q},\ast,\infty\}$ by \eqref{eq:MB} one obtains $M_{f(a'),\ovl{q}} = \gamma_\ell^{2a'}M_{a',\infty} = y\gamma_\ell^{2a'}$. Combining this with application of the first equation of \eqref{eq:MC} yields
\be\label{eq:Mf}  M_{f(a')+t,\ovl{q}+tq} = \gamma_\ell^{-2t} M_{f(a'),\ovl{q}} = y\gamma_\ell^{2a'- 2t} \ .\ee
Letting now $a=f(a')+t$ and $b=\ovl{q}+tq$ and evaluating $t$ and $a'$ in terms of $a$ and $b$ gives $t=bq+\ovl{q}$ and $a'= f^{-1}(a-bq+\ovl{q})$ and hence \eqref{eq:Mf} with $\gamma_\ell^{2q}= \gamma_\ell^{-2}$ reduces to
\be\label{eq:Mff} M_{a,b} = y\gamma_\ell^{2f^{-1}(a-bq+\ovl{q})+ 2b}\ \ {\rm for} \ \ a\ne bq\ ; \ee
the condition $a\ne bq$ is a consequence of $a'\ne \ovl{q}$. Now, \eqref{eq:Mff} together with the information about zero entries of $M$ and about rows and columns marked $\ast$ and $\infty$ determine the matrix $M=M_\ell$ completely. Conversely, one may check that the entries of $M=M_\ell$ that have just been determined satisfy the equations \eqref{eq:MB} and \eqref{eq:MC} with $(-1)^{\delta(\ell)}=(-1)^s$.
\smallskip

With this in hand we are in position to calculate the remaining values of $\chi_\ell$ for $\ell\in L^-$ at the elements $g_4(j)$, represented by the matrix $\Phi_\ell^G(g_4(j))=\off(D_{\rm upp}(j),D_{\rm low}(j))$ given in part (d) of Proposition \ref{prop:miss-rep}; in particular, entries of $D_{\rm upp}(j)$ are zero except for those indexed $(\ast,\infty)$, $(\infty,\ast)$, and $(a,-a-j)$ for $a\in {\rm Ind}_0$, equal, respectively, to $(-1)^\ell\gamma_\ell^{j-1}$, $(-1)^\ell\gamma_\ell^{-j-1}$, and $(-1)^\ell \gamma_\ell^{-j-2a-1}$. By Proposition \ref{prop:eigM} we have $\chi_\ell(g_4(j))={\rm tr}(c^*D_{\rm upp}(j)M)$, and a straightforward evaluation of the trace gives
\[
\chi_\ell(g_4(j)) = i^s\left( \gamma_\ell^{-j} \sum_{a\in {\rm Ind}_0} \gamma_\ell^{-2a} M_{-a-j,a} + (\gamma_\ell^j + \gamma_\ell^{-j})M_{\ast,\infty} \right)\ .
\]
With the help of \eqref{eq:Mff} and $y=\pm q^{-1}$, and also by changing the summation variable from $a$ to $t=-a+(q-1)/2$, the trace equation transforms to
\be\label{eq:trDM}
\chi_\ell(g_4(j)) =\pm i^s q^{-1}\left( \gamma_\ell^{-j} \sum_{t\in {\rm Ind}_0}
\gamma_\ell^{2f^{-1}(t(q+1)-j)} + (\gamma_\ell^j + \gamma_\ell^{-j}) \right)\ .
\ee
Again the key is to use the fact that now $\gamma_\ell^{2(q+1)}=1$, so that to understand the sum in \eqref{eq:trDM} one just needs to study the residues of $f^{-1}(t(q+1)-j)$ mod $(q+1)$, and it is sufficient to do this only for $t\in\{0,1,\ldots,q-2\}$ as the values of $f^{-1}$ depend only on the residue of $t$ but this time mod $(q-1)$. The following auxiliary result will help sort the situation.

\bl\label{lem:f} For $t\in \{0,1,\ldots,q{-}2\}$ the values of $f^{-1}(t(q{+}1)-j)$ are in distinct congruence classes {\rm mod} $(q+1)$; moreover, $f^{-1}(t(q{+}1)-j)\not\equiv 0$ and $\not\equiv j$ {\rm mod} $(q+1)$.\el

\pr Let $e(t)=f^{-1}(t(q+1)-j)$. Suppose that for $t_1,t_2\in\{0,1,\ldots,q-2\}$ one has $e(t_1)=u_1(q+1)+v$ and $e(t_2)=u_2(q+1) +v$ for the same residue $v$ mod $(q+1)$. Observe that the power $z=\xi^{q+1}$ of the primitive element $\xi\in F\cong GF(q^2)$ that has been used throughout is an element of the subfield $F_0\cong FG(q)$ of $F$. By the defining equation for $f$ in part (b) of Proposition \ref{prop:miss-rep} one then has $-1= \xi^{-e(t_1)}+ \xi^{t_1(q+1)-j} = z^{u_1}\xi^v + z^{t_1}\xi^{-j}$ and, similarly, $-1= z^{u_2}\xi^v + z^{t_2}\xi^{-j}$. Eliminating $\xi^v$ from the two equations gives $\xi^{-j}(z^{t_1-u_1} -z^{t_2-u_2}) =z^{-u_2}-z^{-u_1}$. But as $\xi^{-j}\notin F_0$ for $1\le j\le q$ while all the powers of $z$ are in $F_0$, it follows that $u_1 = u_2$ and hence also $t_1=t_2$, proving the first part of the result. The second part is proved similarly, letting $e(t)=j$ and then $e(t)=0$ in the first part of the calculation above. \hfill $\Box$
\medskip

The set ${\rm Ind}_0=\{0,1,\ldots,q^2-2\}$ contains, for each $u\in \{0,1,\ldots,q-2\}$, exactly $q+1$ distinct elements with residue mod $(q-1)$ equal to $u$, and each of these $q+1$ elements of the form $t=u'(q-1)+u$ for $u'\in \{0,1,\ldots, q\}$ determine the same value of the argument of $f^{-1}$ in \eqref{eq:trDM}. Therefore the sum in \eqref{eq:trDM} evaluates to
\be\label{eq:sum-f}  \sum_{t\in {\rm Ind}_0} \gamma_\ell^{2f^{-1}(t(q+1)-j)} = (q+1) \sum_{0\le u\le q-2} \gamma_\ell^{2f^{-1}(u(q+1)-j)}\ . \ee
By Lemma \ref{lem:f}, the values of $f^{-1}$ appearing in \eqref{eq:sum-f} cover exactly $q-1$ out of the possible $q+1$ residue classes mod $(q+1)$, except for the classes $0$ and $j$. But the sum of powers of $\gamma_\ell^2$ with exponents ranging over {\sl all} the $q+1$ distinct residue classes mod $(q+1)$ is equal to zero, because it is precisely the sum of all the $q+1$ distinct complex $(q+1)^{\rm th}$ roots of unity. It follows that the sum of \eqref{eq:sum-f} is equal to the negative of the `missing' two terms corresponding to residue classes $j$ and $0$, that is,
\be\label{eq:sum-ff} \sum_{0\le u\le q-2}\gamma_\ell^{2f^{-1}(u(q+1)-j)}=-\gamma_\ell^{2j}-\gamma_\ell^0 \ . \ee
Substituting now \eqref{eq:sum-ff} into \eqref{eq:sum-f} and then the resulting sum into \eqref{eq:trDM} finally gives the values of $\chi_\ell(g_4(j))$ equal to $\pm i^s(\gamma_\ell^j+\gamma_\ell^{-j})$. We state the corresponding result as a corollary, using the notation $\tau^-_\ell(j) = i^s(\beta_s^j+\beta_s^{-j})$ for $\ell=s(q-1)/2$, $1\le s\le (q-1)/2$, where $\beta_s= \exp(i\pi s/(q{+}1)) = \gamma_\ell$; we will also take into account the earlier observation about zero values of $\chi_\ell(g_3(j))$ for $\ell\in L^-$:

\bc\label{cor:L-} For $\ell=s(q-1)/2\in L^-$, $1\le s\le (q-1)/2$, the pair of characters $\chi_\ell$ and $\chi_\ell\lambda$ are determined by $\chi_\ell(g_3(j)) = 0$ for odd positive $j\le q-2$, and $\chi_\ell(g_4(j)) = \tau^-_\ell(j)$ for odd positive $j\le q$. \hfill $\Box$
\ec

\section{The character table of $M(q^2)$}\label{sec:concl}

It remains to translate the facts derived in Propositions \ref{prop:lift} and Corollaries \ref{cor:L+} and \ref{cor:L-} into a tabular form. In the interest of saving space and displaying all the irreducible characters of $M(q^2)$ in form of a single table one needs to use a number of abbreviations. The following is a legend for reading the Table \ref{tab:Mq2} below.
\smallskip

{\bf Characters.} We use the symbols $\iota$, $\rho$ and $\rho'$ in the first column of Table \ref{tab:Mq2} in their meaning as in Table \ref{tab:PSL} of characters of $\PSL(2,q^2)$, and {\scalebox{0.6}{$\Sigma$}} for the Steinberg character. In the explanations below we will also use the characters $\rho_\ell$ and $\pi_m$ of $\PSL(2,q^2)$. Continuing in the description of character designations in the first column of Table \ref{tab:Mq2}, we let $\rho_{\ell; q} = \rho_\ell + \rho_{\ell q}$ for $\ell\in {\cal U}$ and $\pi_{m;q} = \pi_m + \pi_{mq}$ for $m\in {\cal V}$; the sets ${\cal U}$ with $|{\cal U}|= \frac{1}{8}(q^2-4q+3)$ and ${\cal V}$ with $|{\cal V}|= \frac{1}{8}(q^2-1)$ have been introduced in Section \ref{sec:conj} before Table \ref{tab:untwisted}. The last two entries in the first column of Table \ref{tab:Mq2} are based on the characters $\chi_\ell$ from Corollaries \ref{cor:L+} and \ref{cor:L-}, and we let $\chi^+_\ell=\chi^{\,}_\ell$ for $\ell \in L^+= \{r(q+1)/2;\ 1\le r\le (q-3)/2\}$ and $\chi^-_\ell=\chi^{\,}_\ell$ for $\ell \in L^-= \{s(q-1)/2;\ 1\le s\le (q-1)/2\}$. In rows marked $\psi,\psi\lambda$ for $\psi\in \{\iota,\,${\scalebox{0.6}{$\Sigma$}},\,$\chi^+_\ell, \chi^-_\ell\}$  where double signs appear on some values, those with the bottom signs are assumed to be values of the product of $\psi$ with the alternating character $\lambda$. The number of characters totals to $2+2+1+\frac18(q^2-4q+3) + \frac18(q^2+1) +2\times\frac12(q-3) +2\times\frac12(q-1) = (q+1)(q+5)/4$, which is the number of conjugacy classes as stated at the beginning of section \ref{sec:pre-char}.
\smallskip

{\bf Conjugacy classes.} Again, in the interest of saving space, in displaying representatives of conjugacy classes of $M(q^2)$ in the first row of Table \ref{tab:Mq2} we will omit square brackets, and the second coordinate, $0$ or $1$, will appear as a subscript; thus, for example, instead of $[a(\zeta^j),0]$ and $[\off(\xi^j),1]$ we simply write $a(\zeta^j)_0$ and $\off(\xi^j)_1$. The second and the third row of Table \ref{tab:Mq2} display the number of conjugacy classes ($\#$ cl) and orders of the corresponding centralizers ($|{\rm Cent}|$); the size of a conjugacy class can be obtained by dividing the order $q^2(q^4-1)$ of $M(q^2)$ by the order of the centralizer. The column $a(\xi^j)_0\,({\cal U})$ corresponds to those representatives $a(\xi^j)_0$ for which $j\in {\cal U}$. The column headed by $a(\xi^{j\,'})_0$ actually comprises four columns, namely, those corresponding to the total of $q-2$ values of $j\,'=j(q\pm 1)$ and $j\,'=j(q\pm 1)/2$ from Table \ref{tab:untwisted}, including the bounds and restrictions on $j$ as given in this table. (Observe that although entries in these columns have the same shape, they take {\em different} values of $j$ and $j\,'$ as input.) The column headed $a(\zeta^k)_0$ applies for $k\in {\cal V}$, and the in the last two columns corresponding to $\dia(\xi^j)_1$ and $\off(\xi^j)_1$ the variable $j$ is odd and bounded by $1\le j\le q-2$ and $1\le j\le q$, respectively.
\smallskip

{\bf Entries.} We have used the following in Table \ref{tab:Mq2}: $\omega_\rho(j\ell) = \alpha^{j\ell}+\alpha^{-j\ell} + \alpha^{j\ell q}+ \alpha^{-j\ell q}$, where $\alpha = \exp(4\pi i/(q^2{-}1))$; $\omega_\pi(km) = \beta^{km} + \beta^{-km} + \beta^{kmq} + \beta^{-kmq}$, where $\beta = \exp(4\pi i/(q^2{+}1))$; $\tau^+_\ell(j) = i^r(\alpha_r^j+\alpha_r^{-j})$ for $\ell=r(q+1)/2$, $1\le r\le (q-3)/2$, where $\alpha_r= \exp(i\pi r/(q{-}1))$; and finally $\tau^-_\ell(j) = i^s(\beta_s^j+\beta_s^{-j})$ for $\ell=s(q-1)/2$, $1\le s\le (q-1)/2$, where $\beta_s= \exp(i\pi s/(q{+}1))$. (Note that the complex numbers $\alpha$ and $\beta$ have also been used in Table \ref{tab:PSL}, the character table of $\PSL(2,q^2)$.)
\smallskip

In the terms introduced above we are finally in position to present the character table of the twisted linear group $M(q^2)$.
\medskip

\bt\label{thm:main} The character table of the group $M(q^2)$ is as follows:
\medskip

\begin{small}
\begin{table}[ht]
	\centering
\begin{tabular}{|c||c|c|c|c|c|c|c|c|}
\hline \xrowht[()]{6pt}
$M(q^2)$ & $I_0$ & $u_0$ & $w_0$ & $a(\xi^j)_0\,({\cal U})$ & $a(\xi^{j\,'})_0$ & $a(\zeta^k)_0$ & $\dia(\xi^j)_1$ & $\off(\xi^j)_1$ \\ \hline\hline \xrowht[()]{4pt}
$\#$ cl & $1$ & $1$ & $1$ & {\small $\frac{1}{8}(q^2{-}4q{+}3)$} & $q{-}2$ & $\frac{1}{8}(q^2{-}1)$ & $\frac{1}{2}(q{-}1)$ & $\frac{1}{2}(q{+}1)$ \\ \hline \xrowht[()]{4pt}
$|{\rm Cent}|$ & $q^2(q^4{-}1)$ & $q^2$ & $2(q^2{-}1)$ & $\frac{1}{2}(q^2{-}1)$ & $q^2{-}1$ & $\frac{1}{2}(q^2{+}1)$ & $2(q{-}1)$ & $2(q{+}1)$ \\ \hline\hline \xrowht[()]{4pt}
$\iota,\,\lambda$ & $1$ & $1$ & $1$ & $1$ & $1$ & $1$ & $\pm 1$ & $\pm 1$  \\ \hline \xrowht[()]{4pt}
{\scalebox{0.6}{$\Sigma$}},\,{\scalebox{0.6}{$\Sigma$}}$\lambda$ & $q^2$ & $0$ & $1$ & $1$ & $1$ & $-1$ & $\pm 1$ & $\mp 1$ \\ \hline \xrowht[()]{4pt}
$\rho{+}\rho'$ & $q^2{+}1$ & $1$ & $2$ & $2(-1)^j$ & $2(-1)^j$ & $0$ & $0$ & $0$ \\ \hline \xrowht[()]{4pt}
$\pi_{m;q}$ & $2(q^2{-}1)$ & $-2$ & $0$ & $0$ & $0$ & ${-}\omega_\pi(km)$ & $0$ & $0$\\ \hline \xrowht[()]{4pt}
$\rho_{\ell;q}$ & $2(q^2{+}1)$ & $2$ & $4(-1)^{\ell}$ & $\omega_\rho(j\ell)$ & $\omega_\rho(j\ell)$ & $0$ & $0$ & $0$ \\ \hline \xrowht[()]{4pt}
$\chi^+_\ell,\,\chi^+_\ell\lambda$ & $q^2{+}1$ & $1$ & $2(-1)^{\ell}$ & $\frac{1}{2}\omega_\rho(j\ell)$ & $\frac{1}{2}\omega_\rho(j\ell)$ & $0$ & $\pm\tau^+_\ell(j)$ & $0$ \\ \hline \xrowht[()]{4pt}
$\chi^-_\ell,\,\chi^-_\ell\lambda$ & $q^2{+}1$ & $1$ & $2(-1)^{\ell}$ & $\frac{1}{2}\omega_\rho(j\ell)$ & $\frac{1}{2}\omega_\rho(j\ell)$ & $0$ & $0$ & $\pm\tau^-_\ell(j)$ \\ \hline \hline
\end{tabular}
\caption{Table of irreducible characters of $G=M(q^2)$.}\label{tab:Mq2}
\end{table}
\end{small}

\et

\bigskip

\noindent{\bf Acknowledgements.} The authors thank Gareth A. Jones for helpful discussions and Grahame Erskine for computational results that enabled to verify validity of the character table of $M(q^2)$ for $q\le 13$.
\smallskip

Both authors gratefully acknowledge support from the APVV Research Grants 17-0428 and 19-0308, as well as from the VEGA Research Grants 1/0238/19 and 1/0206/20.

\end{document}